\documentclass[11pt]{article}
\usepackage{mathrsfs}
\usepackage{amsfonts}
\usepackage{amssymb,amsmath,amsthm}
\usepackage{graphicx}
\usepackage{subfigure}
\usepackage{color}
\usepackage{float}

\def\disp{\displaystyle}

\def\dref#1{(\ref{#1})}
\def\crr{\cr\noalign{\vskip2mm}}

\numberwithin{equation}{section}

\newtheorem{theorem}{Theorem}[section]
\newtheorem{definition}[theorem]{Definition}
\newtheorem{lemma}[theorem]{Lemma}
\newtheorem{remark}[theorem]{Remark}

\newcommand{\mm}     {{\hbox{\hskip 0.5pt}}}

\newcommand{\bluff}  {{\hbox{\raise 15pt \hbox{\mm}}}}
\newcommand{\sbluff} {{\hbox{\raise  9pt \hbox{\mm}}}}

\newcommand{\sbm}[1]{\left[\begin{smallmatrix} #1
   \end{smallmatrix}\right]}

\begin{document}
\renewcommand{\thefootnote}{\fnsymbol{footnote}}
\renewcommand{\thefootnote}{\fnsymbol{footnote}}
\newcommand{\footremember}[2]{%
   \footnote{#2}
    \newcounter{#1}
    \setcounter{#1}{\value{footnote}}%
}
\newcommand{\footrecall}[1]{%
    \footnotemark[\value{#1}]%
}
\makeatletter
\def\blfootnote{\gdef\@thefnmark{}\@footnotetext}
\makeatother
\begin{center}
{\LARGE \bf Disturbance Observer-Based Boundary Control   \\ 
for an Anti-Stable Stochastic Heat Equation with  \\[0.6ex]
  Unknown Disturbance
}\\[4ex]
Ze-Hao Wu, Hua-Cheng~Zhou, Feiqi Deng, and Bao-Zhu~Guo 
\blfootnote{This work was supported by the National Natural Science Foundation of China.}
\blfootnote{Z.H. Wu (zehaowu@amss.ac.cn) is with School of Mathematics and Big Data, Foshan University, Foshan 528000;
H.C. Zhou (hczhou@amss.ac.cn) is with the School of Mathematics and Statistics,
 Central South University, Changsha 410075; 
Feiqi Deng is with Systems Engineering Institute, South China University of Technology, Guangzhou 510640;
B.Z. Guo (bzguo@iss.ac.cn) is with Key Laboratory of System and Control,
        Academy of Mathematics and Systems Science,
        Chinese Academy of Sciences, Beijing 100190, China. 
 }
\end{center}
\vspace{3mm}

\noindent {\bf Abstract:} In this paper, a novel control strategy namely  disturbance observer-based control
is first applied to     stabilization and  disturbance
rejection for an anti-stable stochastic heat
 equation with Neumann boundary actuation and unknown boundary external disturbance generated by an exogenous system.
A disturbance observer-based boundary control is designed based on the backstepping approach and estimation/cancellation
strategy, where the unknown disturbance is estimated in real time
by a disturbance observer  and { rejected in the closed-loop, while the in-domain multiplicative noise
whose intensity is within a known finite interval is attenuated}. It is shown  that
 the resulting closed-loop system is    exponentially stable in the sense of both mean square and
 almost surely. A numerical example is demonstrated  to validate the effectiveness of the proposed control approach. \vspace{3mm}

\noindent {\bf Keywords:} Stochastic heat equation; boundary control; disturbance rejection;  stabilization;
 backstepping approach.
\vspace{3mm}

\noindent {\bf AMS subject classifications:} 35K05, 37L55, 93B52.

\section{Introduction} 
{D}{isturbances} are  ubiquitous in  many practical control
systems, which often cause   negative effects  on
 performance of the control plant.
 For the sake of control precision, many control approaches have been
developed since 1970s to cope with disturbances
in term of disturbance attenuation or disturbance rejection.
The stochastic control 
and robust control 
  are two of representative disturbance
attenuation approaches, where the former is often used
for attenuating noises
with known statistical characteristics while the latter can deal with
 more general disturbances.
However, most of the  robust control approaches are  on the
worst case scenario, which may lead to control design rather conservative.
Based on estimation/cancellation strategy, some novel active anti-disturbance control approaches
like disturbance observer based control (DOBC) \cite{whchen2004,li2016}
have been proposed for disturbance rejection for control systems over the past two decades.
The core idea of these active anti-disturbance control approaches is that the disturbances affecting system performance
 can be estimated by a disturbance observer and then be  compensated in the closed-loop.
Owing to the estimation/cancellation characteristics, the active anti-disturbance control is
capable of eliminating
the disturbances before negative effects  are caused and  at the same time,  the control energy can
be reduced significantly in engineering applications.

The disturbance rejection for distributed parameter systems by active anti-disturbance boundary control approaches
has been paid increasing attention in the last  two decades, see, for instance \cite{wuhn20162,ren2018,zhoutac2019,Deu2011,Deutscher2015auto} and the references therein.
The DOBC approach to the  stabilization  of
nonlinear parabolic partial differential equation (PDE) systems subject to external disturbance has been investigated in \cite{wuhn20162}.
The uncertainty and disturbance estimator (UDE)-based robust control approach to the stabilization of an unstable parabolic PDE
 with a Dirichlet type boundary actuator and an unknown time-varying input disturbance has been addressed in \cite{ren2018}. 
An infinite-dimensional observer-based output feedback boundary control has been designed for 
a multi-dimensional heat equation subject to
boundary unmatched disturbance in  \cite{zhoutac2019}.
{The output regulation has been developed for linear distributed-parameter systems by finite-dimensional
dual observers \cite{Deu2011} and parabolic PDEs by a backstepping approach
\cite{Deutscher2015auto}, respectively.}

Nevertheless, the external  disturbance appears  most often in  random way
in practice, which is neglected in literatures aforementioned.
Actually, there have been many control
designs for  finite-dimensional stochastic systems driven by white noise,
see, for instance \cite{pan1999,krstic1} and the references therein.
Specially, some active anti-disturbance control methods to
 disturbance rejection for finite-dimensional stochastic systems have been proposed.
For example, problems of the composite DOBC and $H_{\infty}$ control  for Markovian jump systems
and the DOBC  for a class of stochastic systems with multiple disturbances
have been studied in \cite{dobcs1} and \cite{dobcs2}, respectively; An extended state observer-based  output feedback
stabilizing control has been designed for
 a class of stochastic systems subject to bounded stochastic noise \cite{Gwz2016}.

As  one of the active  anti-disturbance control approaches,
the DOBC has been widely applied in engineering applications with good disturbance rejection
performance and robustness, see, for instance \cite{appliation1,appliation2} and the references therein.
 However, there is still no relevant study from theoretical perspective on DOBC for stochastic distributed parameter systems.
In this paper,  we demonstrate for the first time, through  an
anti-stable stochastic heat equation with unknown boundary external disturbance, the DOBC approach
to  stabilization and disturbance rejection for stochastic distributed parameter systems.
The main contributions and novelty of this paper can be summarized  as follows. a)~
From a theoretical perspective, the applicability of the powerful DOBC control technology is first
expanded to a class of stochastic distributed parameter systems with unknown boundary external disturbance;
b)~The unknown boundary external disturbance   is rejected completely by virtue of estimation/cancellation strategy of the DOBC approach,
{while the in-domain multiplicative noise with bounded  intensity is attenuated};
c)~Not only the mean square exponential stability but also
  the almost surely exponential stability are obtained for the resulting closed-loop system.

We proceed as follows. In the next section, section \ref{Se2}, some problem formulation and preliminaries
are presented. In section \ref{Se3}, both  design of the DOBC boundary control and  stability
 of the closed-loop system are  discussed and stated. A numerical example  is
presented in section \ref{Se4}, followed up concluding remarks in  section \ref{Se5}. The proofs of the
main results are arranged in Appendix.

\section{Problem formulation and preliminaries}\label{Se2}

We first introduce some notations.  The $I_{n}$ denotes the $n$-dimensional identity matrix and
$L^{2}(0,1)$ is the space of all real-valued functions that
are square Lebesgue integrable over  $(0,1)$.
Let $(\Omega,\mathcal{F},\mathbb{F}, P)$ be
a complete filtered probability space with a filtration
$\mathbb{F}=\{\mathcal{F}_{t}\}_{t\geq 0}$  on which
a one-dimensional standard Brownian motion $B(t)$ is defined.
{A stochastic process $f(t,\omega):[0,\infty)\times \Omega\rightarrow \mathbb{R}^{n}$
is called $\mathbb{F}$-adapted if for every  $t\geq 0$, the function $\omega\rightarrow f(t,\omega)$
is $\mathcal{F}_{t}$-measurable. For notation simplicity we use  $f(t)$ to denote a stochastic process $f(t,\omega)$.}
Let $V$ be a Banach space. A
sub-$\sigma$-algebra $\mathcal{M}$ of $\mathcal{F}$, denoted
by $L^{2}_{\mathcal{M}}(\Omega;N)$, is the set of all $\mathcal{M}$-measurable ($N$-valued)
random variables $f:\Omega \rightarrow N$ such that $\mathbb{E}|f|^{2}_{N}<\infty$.
Let $ H^1(0,1)$ and $H^2(0,1)$ be the Sobolev spaces. Set
$
L^2_{\mathbb{F}}(0,T;L^2(\Omega; V))=\{f:(0,T)\times\Omega\rightarrow V| f(\cdot)$ is $\mathbb{F}\mbox{-adapted
and}$ 
$\int^{T}_{0}(\mathbb{E}|f(t)|^{2}_{V})dt<\infty\}$,
$C_{\mathbb{F}}([0,T];L^2(\Omega;V))$ $=\{f:[0,T]\times\Omega\rightarrow V| f(\cdot)$ is 
$\mathbb{F}\mbox{-adapted
and}$
$(\mathbb{E}|f(t)|^{2}_{V})^{\frac{1}{2}}$
is continuous\},
$C_{\mathbb{F}}(0,\infty;$
$L^2(\Omega; V))$
$=\{f:(0,\infty)\times\Omega\rightarrow V| f(\cdot)$\; $\mbox{is}\;\mathbb{F}$-adapted and 
$(\mathbb{E}|f(t)|^{2}_{V})^{\frac{1}{2}}\; \mbox{is continuous}\}$,\;
$C^1_{\mathbb{F}}(0,T;L^2(\Omega;V))=\{f:(0,T)\times\Omega\rightarrow V| f(\cdot)\; \mbox{is}\;\mathbb{F}\mbox{-adapted
and}(\mathbb{E}|f(t)|^{2}_{V})^{\frac{1}{2}}$ 
is continuously differentiable\}. 
All the above spaces are endowed with the usual canonical norms.

In this paper, we consider  stabilization and disturbance rejection  for a one-dimensional anti-stable stochastic heat equation
driven by multiplicative white noise with unknown  boundary external disturbance  as follows:
\begin{equation}\label{SPDE-Heat}
 \left\{\begin{array}{l}
dy(x,t)=y_{xx}(x,t)dt+a(x)y(x,t)dt+\sigma y(x,t)dB(t), \cr
y_x(0,t)=0,\; t\geq0,\cr y_x(1,t)=u(t)+w(t),\; t\geq0,\cr
y(x,0)=y_0(x),\;0\leq x\leq1,
\end{array}\right.
\end{equation}
where $y(x,t)$ is the {system state} representing the temperature
profile at the spatial position $x\in [0,1]$ and the time $t\in [0,\infty)$,
$a(\cdot)\in L^2(0,1)$,
$\sigma$ is a constant representing the intensity of the multiplicative white noise {with a known upper bound for its absolute value},
 $y_0(\cdot)\in L^2_{\mathcal{F}_0}(\Omega;L^2(0,1))$ is the initial value,
and $u(t)$ is the boundary control input,   
 $w(t)\in \mathbb{R}$ is the unknown disturbance which could be the temperature perturbation generated from an exogenous system as follows:
\begin{equation}\label{exogenous}
 \left\{\begin{array}{l}
\dot{\xi}(t)=A\xi(t), \crr\disp
w(t)=C\xi(t),
\end{array}\right.
\end{equation}
where $\xi(t)\in \mathbb{R}^{n}$ is an unknown exogenous signal, and $A\in \mathbb{R}^{n\times n}$, $C\in \mathbb{R}^{1\times n}$ are known matrices.
{Throughout the paper, the stochastic differentials including $dy(x,t)$ are with respect to the time  $t$. }

The following Assumption (A1) is for  estimation  of the unknown disturbance  $w(t)$.
\hspace{0.1cm}

{\bf Assumption (A1).} For  the exogenous  system (\ref{exogenous}), the  pair $(A, C)$ is supposed to be observable.

\begin{remark}
 The exogenous system  (\ref{exogenous}) is a typical form for disturbance in
 output regulation, which covers all finite harmonic disturbance with
 known frequency but unknown amplitude and
phase.   
The harmonic disturbance can be considered as an approximation of
periodic disturbance (\cite{whchen2004,dobcs1,application2}) and has been discussed intensively   by means of  
 the DOBC approach (see, e.g., \cite[p.47]{li2016}).
\end{remark}

{Similar to both mean square exponential stability and almost surely exponential stability of stochastic differential equations (see, e.g.,\cite{mao2007}), we introduce the following stability for   system \dref{SPDE-Heat}.
\begin{definition}
System \dref{SPDE-Heat} is said to be  mean square exponentially stable if there are positive constants $M$ and $\eta$
such that
$$
\mathbb{E}|y(\cdot,t)|^{2}_{L^2(0,1)}\leq Me^{-\eta t},\; \forall t\geq 0,
$$
for all initial value  $y_0\in L^2_{\mathcal{F}_0}(\Omega;L^2(0,1))$, where $M$ depends on the initial value $y_0$.
\end{definition}

\begin{definition}
System \dref{SPDE-Heat} is said to be almost surely exponentially stable if
$$
\limsup_{t\to \infty}\frac{1}{t}\log|y(\cdot,t)|_{L^2(0,1)}<0 \;\;\;\; \mbox{almost surely,}
$$
for all initial value $y_0\in L^2_{\mathcal{F}_0}(\Omega;L^2(0,1))$.
\end{definition}
}

\begin{remark}\label{Rem-exam}
System \dref{SPDE-Heat} without boundary control could be neither stable in mean
square nor stable in almost surely even if the  boundary input is vanishing: $u(t)=w(t)\equiv0$. Actually, by taking $a(x)=4\pi^2+1.005$, $\sigma=0.1$ and $y_0(x)=\cos(2\pi x)$
in \dref{SPDE-Heat}, $y(x,t)=\cos(2\pi x)e^{t+0.1B(t)}$ solves \dref{SPDE-Heat}. Since $\lim_{t\to\infty}\frac{B(t)}{t}=0$ almost surely,  it is clear that
$\int_0^1 y^2(x,t)dx\rightarrow \infty$ almost surely as $t\rightarrow\infty$.
Moreover, $\mathbb{E}\int_0^1y^2(x,t)dx=\int_0^1y^2_0(x)dx\cdot \mathbb{E}e^{2t+0.2B(t)}
=\int_0^1y^2_0(x)dx\cdot e^{2.02t}\rightarrow \infty$ as $t\rightarrow\infty$.
\end{remark}

In what follows, we consider  system \dref{SPDE-Heat} in the state space $L^2_{\mathbb{F}}(\Omega; L^2(0,1))$  with the usual canonical norm.
The control objective is to design a  disturbance observer-based boundary control so that
the closed-loop system is exponentially stable in both mean square and almost surely.

\section{DOBC boundary control design and main results}\label{Se3}

\begin{figure}[ht]\centering
 {\includegraphics[width=14.5cm,height=5cm]{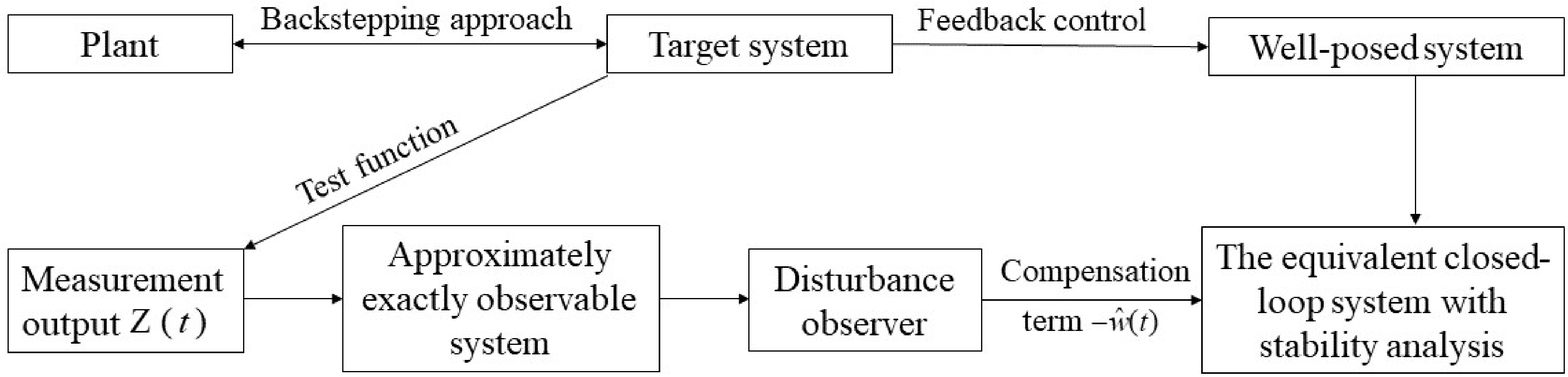}\label{Fig-openy-b--}}
\caption{The framework of the DOBC boundary control design and theoretical approach.}
\label{Fig-ref-cblock}
\end{figure}

The framework of the DOBC boundary control design and theoretical approach
can be simply explained  in  Figure \ref{Fig-ref-cblock}.
By the backstepping approach, the controlled system
is first  {transformed} into an equivalent  ``good" target system for which we are able to design
a feedback stabilizing control without unknown boundary external  disturbance being considered.
Second, by a test function, a measurement output is introduced which is the
solution of an approximatively exactly observable It\^{o}-type stochastic system with $w(t)$ being its
external disturbance, and then a finite-dimensional disturbance observer is designed for real-time estimation of
 the unknown disturbance $w(t)$. The DOBC boundary control constructed by the feedback
stabilizing control  and a compensation term
by the estimate of the unknown disturbance is finally  designed to obtain the resulting closed-loop system.

Motivated by \cite{AS-Krstic-book},
we introduce an invertible transformation $\Lambda: y\in L^2_{\mathbb{F}}(\Omega;L^2(0,1))\to z\in L^2_{\mathbb{F}}(\Omega;L^2(0,1))$ as follows:
\begin{equation}\label{bstep}
z(x,t)=y(x,t)-\int_0^xk(x,\zeta)y(\zeta,t)d\zeta\triangleq\Lambda y(x,t),
\end{equation}
where {$k(\cdot,\cdot)$ defined on $G:=\{(x,\zeta)\in \mathbb{R}^{2}: 0\leq \zeta\leq x\leq 1\}$} is the kernel function specified in  (\ref{kerequ}) later.
Let $c$  be a positive constant. By It\^{o}'s differentiation rule, a direct computation shows that
\begin{eqnarray}\label{induker}
\hspace{-0.5cm}&&dz(x,t)-z_{xx}(x,t)dt+cz(x,t)dt-\sigma z(x,t)dB(t)\cr \hspace{-0.5cm}&&
 =dy(x,t)-\big[k(x,x)y_{x}(x,t)-k(x,0)y_{x}(0,t)\cr\hspace{-0.5cm} &&
 -(k_{\zeta}(x,x)y(x,t)-k_{\zeta}(x,0)y(0,t))\big]dt\cr \hspace{-0.5cm}&&
 -\int_0^xk_{\zeta\zeta}(x,\zeta)y(\zeta,t)d\zeta dt-\int_0^xk(x,\zeta)a(\zeta)y(\zeta,t)d\zeta dt\cr\hspace{-0.5cm}&&-\sigma\int_0^xk(x,\zeta)y(\zeta,t)d\zeta dB(t)\cr \hspace{-0.5cm}&&
 -\big[y_{xx}(x,t)-\frac{dk(x,x)}{dx}y(x,t)-k(x,x)y_x(x,t)\big]dt\cr\hspace{-0.5cm} &&
 -\big[-k_{x}(x,x)y(x,t)-\int_0^xk_{xx}(x,\zeta)y(\zeta,t)d\zeta\big]dt\cr \hspace{-0.5cm}&&
  +c\big[y(x,t)-\int_0^xk(x,\zeta)y(\zeta,t)d\zeta\big]dt-\sigma
z(x,t)dB(t)\cr \hspace{-0.5cm}&&
=\int_0^x[k_{xx}(x,\zeta)-k_{\zeta\zeta}(x,\zeta)-(c+a(\zeta))k(x,\zeta)]y(\zeta,t)d\zeta dt\cr\hspace{-0.5cm} &&
+\big[\frac{dk(x,x)}{dx}+k_{x}(x,x)+k_{\zeta}(x,x)
+a(x)+c\big]y(x,t)dt\cr \hspace{-0.5cm}&&-k_\zeta(x,0)y(0,t)dt.
\end{eqnarray}

 The kernel function $k(\cdot)$ is chosen to satisfy
  the following partial differential equation (PDE):
\begin{equation}\label{kerequ}
\left\{\begin{array}{l}
k_{xx}(x,\zeta)-k_{\zeta\zeta}(x,\zeta)=(c+a(\zeta))k(x,\zeta),\cr
k_\zeta(x,0)=0,\cr
\frac{dk(x,x)}{dx}+k_{x}(x,x)+k_{\zeta}(x,x)=-a(x)-c.
\end{array}\right.
\end{equation}
It  can be seen that the right side of \dref{induker} becomes zero.
 Set $z_x(0,t)=0$. Then the boundary condition in \dref{SPDE-Heat} yields
$0=z_x(0,t)=y_x(0,t)-k(0,0)y(0,t)=-k(0,0)y(0,t)$,
which implies $k(0,0)=0$.  This, together with the third equation
in \dref{kerequ},
indicates that the  PDE \dref{kerequ}  becomes
\begin{equation}\label{kerequPDE}
\left\{\begin{array}{l}
k_{xx}(x,\zeta)-k_{\zeta\zeta}(x,\zeta)=(c+a(\zeta))k(x,\zeta),\cr
k_\zeta(x,0)=0,\cr
k(x,x)=-\frac{1}{2}\int_0^xa(\zeta)d\zeta-\frac{c}{2}x.
\end{array}\right.
\end{equation}

\begin{remark}
 By  \cite[Theorem 2.1]{AS-Krstic-book},
\dref{kerequPDE} admits a unique solution which is twice
continuously differentiable in the domain  $G= \{(x,\zeta)\in \mathbb{R}^{2}: 0\leq \zeta\leq x\leq 1\}$.
In particular,  if $a(\cdot)$ is a constant function, for example, $a(\zeta)\equiv 0$,  then $k(x,\zeta)$ can be found analytically as
$k(x,\zeta)=-cx\frac{I_1(\sqrt{c(x^2-\zeta^{2})})}{\sqrt{c(x^2-\zeta^{2})}}$,
where $I_1(\cdot)$ is a first-order modified Bessel function given by
$I_1(\varpi)=\sum^{\infty}_{i=0}\frac{\varpi^{2i+1}}{2^{2i+1}i!(i+1)!}.$
\end{remark}

The following Lemma \ref{lemma3} brought from  Theorem 2.2 of  \cite{AS-Krstic-book}
demonstrates that the transformation operator $\Lambda$ defined in \dref{bstep} is invertible.

\begin{lemma}(\cite{AS-Krstic-book})\label{lemma3}
The linear operator $\Lambda:L^2_{\mathbb{F}}(\Omega;L^2(0,1))\rightarrow L^2_{\mathbb{F}}(\Omega;L^2(0,1))$ defined in \dref{bstep} is bounded
invertible, and
\begin{equation}\label{bstep-inv}
y(x,t)=z(x,t)+\int_0^xl(x,\zeta)z(\zeta,t)d\zeta,
\end{equation}
where $l\in C^2(G)$.
\end{lemma}

 Under the  invertible transformation \dref{bstep}, system \dref{SPDE-Heat} is
transformed  into the following equivalent one:
\begin{equation}\label{SPDE-Heat-trans}
\left\{\begin{array}{l}
dz(x,t)=z_{xx}(x,t)dt-cz(x,t)dt+{\sigma}z(x,t)dB(t),\cr
 z_x(0,t)=0, \; t\geq0,\cr
z_x(1,t)=u(t)+w(t)-k(1,1)y(1,t)\cr \hspace{1.4cm}-\int_0^1k_x(1,\zeta)y(\zeta,t)d\zeta, \;
t\geq0,\cr z(x,0)=z_0(x),\;0\leq x\leq1.
\end{array}\right.
\end{equation}

By introducing  a new control law $u_0(t)$, we design
\begin{equation}\label{311equation}
u(t)=u_0(t)+k(1,1)y(1,t)+\int_0^1k_x(1,\zeta)y(\zeta,t)d\zeta,
\end{equation}
under which  system \dref{SPDE-Heat-trans} becomes
 the  ``good" target system we are looking for because its
 well-posedness and stability can be easily obtained for
 $u_{0}(t)+w(t)\equiv 0$. Hence, a key step is to estimate the unknown disturbance $w(t)$  first by a disturbance observer and
eliminate it  in feedback loop  by the DOBC boundary control. In order to estimate the unknown
 disturbance $w(t)$ in real time, we introduce { an average-type signal:}
\begin{equation}\label{zequation}
Z(t)=\int_0^1\cos(\pi x) z(x,t)dx,
\end{equation}
 {where $\cos(\pi x)$ is a test function.} By applying  It\^{o}'s differentiation rule to \dref{zequation},
 the unknown boundary external disturbance $w(t)$
can be {transformed} into the following one-dimensional It\^{o}-type stochastic system:
\begin{equation}\label{Z-sde-U0}
dZ(t)=-\big[(c+\pi^2)Z(t)+u_0(t)+w(t)\big]dt+\sigma Z(t)dB(t).
\end{equation}
We then design a finite-dimensional disturbance observer as follows:
\begin{equation}\label{dobc-part}
\left\{\begin{array}{l}
\dot{\vartheta}(t)=(A+LC)\vartheta(t)+(A+LC)LZ(t)\cr \hspace{1cm}+L[(c+\pi^2)Z(t)+u_0(t)],\cr
\widehat{\xi}(t)=\vartheta(t)+LZ(t),\cr
\widehat{w}(t)=C\widehat{\xi}(t),
\end{array}\right.
\end{equation}
where $\vartheta(t)$ is the internal state of the disturbance observer \dref{dobc-part}, $L\in \mathbb{R}^{n\times 1}$
is chosen so that $A+LC$ is Hurwitz by Assumption (A1), and
$\widehat{\xi}(t)$ and $\widehat{w}(t)$ are the estimates of $\xi(t)$ and the unknown disturbance $w(t)$, respectively.
Set  the observer estimation error as $\eta(t)=\widehat{\xi}(t)-\xi(t)$.
It is easy to verify that $\eta(t)$ satisfies the following It\^{o}-type stochastic differential equation
\begin{equation}\label{e-sde-F}
d\eta(t)=(A+LC)\eta(t)dt+\sigma LZ(t)dB(t).
\end{equation}
{
\begin{remark}It is pointed out that
the proposed finite-dimensional disturbance observer for real-time estimation of unknown
boundary external disturbance   of stochastic PDEs
is systematic. Actually, let the state of the  stochastic PDE also be denoted by $z(x,t)$.
The idea is to find  an average-type output $Z(t)=\int_0^1h(x)z(x,t)dx$ where $h(x)$ is a test function  governed by the stochastic PDE
so that $Z(t)$ satisfies an approximatively exactly observable It\^{o}-type stochastic system with $w(t)$ being the external disturbance
 determined uniquely by the  average-type output $Z(t)$.
Specifically, for the case that $z(x,t)$ satisfies   (\ref{SPDE-Heat-trans}) with $u(t)$ being designed in (\ref{311equation}), it can be
obtained that the average-type output $Z(t)=\int_0^1h(x)z(x,t)dx$ satisfies
\begin{eqnarray}\label{derivedsys}
\begin{array}{l}\disp
dZ(t)=h(1)[u_0(t)+w(t)]dt+[-z(1,t)h'(1)+z(0,t)h'(0)]dt\crr\disp
 \ \ \ \ \ \ \ \ +\int_0^1h''(x)z(x,t)dxdt-cZ(t)dt+\sigma Z(t)dB(t).
\end{array}
\end{eqnarray}
 Hence, to guarantee that  $w(t)$ is uniquely determined by the average-type output $Z(t)$,
 the test function $h(x)$ can be any one that satisfies the condition
\begin{equation}\label{condition}
h'(0)=h'(1)=0, h(1)\neq0, h''(x)=\nu h(x),
\end{equation}
where $\nu$ is any constant. This is the reason we  choose $h(x)=\cos(\pi x)$ for our case. 
In fact, for any $T>0$,
\begin{equation}
Z(t)\equiv 0,u_0(t)\equiv0, \forall t \in [0,T]\Rightarrow w(t)\equiv 0, \forall t \in [0,T],
\end{equation}
which means that $w(t)$ is uniquely determined by the average-type output $Z(t)$. This is
similar to the exact observability of It\^{o}-type stochastic systems without external disturbance that the state is uniquely
determined by output, see for instance \cite{exactobservab}.  Therefore, a natural way is to use the output $Z(t)$ of system (\ref{Z-sde-U0})
 to design a disturbance observer to estimate $w(t)$, and hence the idea of  finite-dimensional disturbance observer
 can be adopted to design the observer (\ref{dobc-part}).
In addition, no matter what  test function satisfying condition (\ref{condition}) to be chosen,
 the design framework of the disturbance observer  does not change for
 an approximatively exactly observable linear It\^{o}-type stochastic system (\ref{derivedsys}).
 The control design also does not change because
in the steady state, $\hat{w}(t)$ approaches $w(t)$ in both mean square and almost sure sense.
\end{remark}
}

Since $A+LC$ is Hurwitz, the mean square exponential stability of the error system \dref{e-sde-F} can be expected
if the intensity of the white noise is ``small" to some extent,
and then the  almost surely exponential stability can also be concluded directly due to the common linear growth condition.
Thus, we can design  the control law $u_0(t)=-\widehat{w}(t)$
 in (\ref{311equation}) to cancel the unknown disturbance $w(t)$ in real time.
We first consider the mean square exponential stability for  $(Z(t),\eta^{\top}(t))^{\top}$ which is governed by
\begin{equation}\label{coupled-Z-e-sys}
d\sbm{Z(t)\\ \eta(t)}=\sbm{-(c+\pi^2)&C\\ 0& A+LC }\sbm{Z(t)\\ \eta(t)}dt+\sbm{\sigma&0\\ \sigma L&0 }\sbm{Z(t)\\ \eta(t)}dB(t).
\end{equation}
{By Assumption (A1),  $\sbm{-(c+\pi^2)&C\\ 0& A+LC }$ is Hurwitz and hence there exists a unique positive matrix $Q_{c}:=\sbm{Q_{c1}&Q_{c2}\\ Q_{c3}& Q_{c4} }\in\mathbb{R}^{(n+1)\times(n+1)}$
depending  on $c$ such that
\begin{equation}\label{Ly-alg-equ}
\sbm{-(c+\pi^2)&C\\ 0& A+LC}^{\top}Q_{c}+Q_{c}\sbm{-(c+\pi^2)&C\\ 0& A+LC}=-\mathbb{I}_{(n+1)},
\end{equation}
where $Q_{c1}\in \mathbb{R}^{+},Q_{c2}=Q^{\top}_{c3}\in \mathbb{R}^{1\times n}$, and $Q_{c4}\in \mathbb{R}^{n\times n}$.
A direct computation can easily show that $Q_{c1}=\frac{1}{2(c+\pi^{2})}$, $Q_{c2}=Q^{\top}_{c3}=-\frac{1}{2(c+\pi^{2})}C[(A+LC)-(c+\pi^{2})\mathbb{I}_{n}]^{-1}$, and $Q_{c4}$ satisfies $(A+LC)^{\top}Q_{c4}+Q_{c4}(A+LC)+C^{\top}Q_{c2}+Q_{c3}C=-\mathbb{I}_{n}$.
It can be further obtained that
\begin{eqnarray}\label{ucontant}
&& \mu_{c}:=\lambda_{\max}\left(\sbm{1 &0\\  L&0 }^{\top}Q_{c}\sbm{1&0\\  L&0 }\right)
\cr&&
=Q_{c1}+Q_{c2}L+L^{\top}Q_{c3}+L^{\top}Q_{c4}L,
\end{eqnarray}
which is a  $c$-dependent  positive constant because $Q_c$ is  a positive definite matrix and $\sbm{1 &0\\  L&0 }^{\top}Q_{c}\sbm{1&0\\  L&0 }$ is a nonzero positive semi-definite one. }
The mean square and almost surely exponential stability of system \dref{coupled-Z-e-sys} can be summarized as the following Lemma \ref{lem-Z-e-to0}.

\begin{lemma}\label{lem-Z-e-to0}
{Suppose that  Assumption (A1)  and $|\sigma|<\frac{1}{\sqrt{\mu_{c}}}$ hold.}
Then, the solution of system \dref{coupled-Z-e-sys} satisfies
\begin{eqnarray}\label{estiomerror}
\mbox{(i)}\hspace{-0.5cm}&&\mathbb{E}[|Z(t)|^2+|\eta(t)|^2]\cr\hspace{-0.5cm} &&\leq\frac{ \lambda_{\max}({Q_{c}})}{\lambda_{\min}({Q_{c}})}\mathbb{E}[|Z(0)|^2+|\eta(0)|^2]e^{-\frac{(1-{\mu_{c}}\sigma^2)}{\lambda_{\max}({Q_{c}})}t},\; \forall t\geq 0;
\end{eqnarray}
\begin{eqnarray}\label{almostlysure}
\hspace{-0.5cm}\mbox{(ii)}\hspace{-0.5cm}&&\max\{\limsup_{t\to \infty}\frac{1}{t}\log|Z(t)|,\;\limsup_{t\to \infty}\frac{1}{t}\log|\eta(t)|\}\cr
\hspace{-0.6cm}&&\leq -\frac{(1-{\mu_{c}}\sigma^2)}{2\lambda_{\max}({Q_{c}})}\;\;\;\;\; \mbox{almost surely.}
\end{eqnarray}
\end{lemma}

{\bf Proof.}  See {\bf ``Proof of Lemma \ref{lem-Z-e-to0}''} in Appendix A.

The DOBC boundary control is
\begin{equation}\label{controllerde}
u(t)=k(1,1)y(1,t)+\int_0^1k_x(1,\zeta)y(\zeta,t)d\zeta-\widehat{w}(t),
\end{equation}
under which system \dref{SPDE-Heat-trans} becomes
\begin{equation}\label{SPDE-Heat-trans-U0-closed}
\left\{\begin{array}{l}
dz(x,t)=z_{xx}(x,t)dt-cz(x,t)dt+{\sigma}z(x,t)dB(t),\cr
 z_x(0,t)=0, \; t\geq0,\cr
z_x(1,t)=-\widehat{w}(t)+w(t)\triangleq \tilde{w}(t),\;
t\geq0,\cr z(x,0)=z_0(x),\;0\leq x\leq1.
\end{array}\right.
\end{equation}
{ It is noted that in order to design disturbance observer (\ref{dobc-part}) and DOBC boundary control \dref{controllerde},
we used actually three measured outputs $Y(t):=\{y(1,t),Z(t),\int_0^1k_x(1,\zeta)y(\zeta,t)d\zeta\}$,
where $y(1,t)$ is a pointwise measurement output, and $Z(t)$ defined in \dref{zequation} and $\int_0^1k_x(1,\zeta)$ $y(\zeta,t)d\zeta$ with $k(\cdot,\cdot)$ specified
in (\ref{kerequ}) are two average-type measurement outputs. Such kinds of measurement outputs can be found in Chapters 4 or 5
of the monograph \cite{book1995Hans}  where they are regarded as an average measurement of the temperature around a
certain point. In a recent paper \cite{Fridman1}, these signals are regarded as the non-local ones.
 In addition, by choosing appropriately
 measurement distribution functions, they can be realized by a finite number of sensors in practice, see for instance \cite{Deu2011,wustochasticpde}.
Nevertheless, we consider present paper as state feedback. The output feedback stabilization is a separate issue
based on state feedback stabilization.}

Next, we first show the well-posedness of the closed-loop system (\ref{SPDE-Heat-trans-U0-closed}).
\begin{lemma}\label{lem-existencewell-u-sys}
{Suppose that  Assumption (A1) and $|\sigma|<\frac{1}{\sqrt{\mu_{c}}}$ hold.}
 Then,  for any initial value $z_0\in L^2_{\mathcal{F}_0}(\Omega;L^2(0,1))$, the closed-loop system \dref{SPDE-Heat-trans-U0-closed} admits a
   unique {weak} solution
$z\in C_{\mathbb{F}}(0,+\infty;L^2(\Omega; L^2(0,1)))$; Moreover,
for any $T>0$, there exists a positive constant $C(T)$ such that
\begin{eqnarray}\label{u-theta-estinq}
 && |z|_{C_{\mathbb{F}}([0,T];L^2(\Omega;
L^2(0,1)))}+|z|_{L^2_{\mathbb{F}}(0,T;L^2(\Omega; L^2(0,1)))}\cr &&  \leq
C(T)\left[|z_0|_{L^2_{\mathcal{F}_0}(\Omega;L^2(0,1))}+|\widetilde{w}|_{ L^2_{\mathbb{F}}(0,T;L^2(\Omega;\mathbb{R}))}\right].
\end{eqnarray}
\end{lemma}
{\bf Proof.}  See {\bf ``Proof of Lemma \ref{lem-existencewell-u-sys}''} in  Appendix B.

{ The positive constant $\theta^{*}$ used hereinbelow is defined as
\begin{equation}\label{canostff}
\hspace{-0.1cm}\theta^{*}=\left\{\begin{array}{l}
\min\{2c-2-3\sigma^{2},\frac{1-\mu_{c}\sigma^{2}}{\lambda_{\max}(Q_{c})}\},\mbox{if}\; 2c-2-3\sigma^{2}\neq\frac{1-\mu_{c}\sigma^{2}}{\lambda_{\max}(Q_{c})},\cr
\theta,\hspace{3.4cm}\mbox{if}\; 2c-2-3\sigma^{2}=\frac{1-\mu_{c}\sigma^{2}}{\lambda_{\max}(Q_{c})},
\end{array}\right.
\end{equation}
where $\theta$ is a fixed positive constant satisfying $\theta<2c-2-3\sigma^{2}$.}

The mean square exponential stability of the closed-loop system  \dref{SPDE-Heat-trans-U0-closed} is
summarized in the following Lemma \ref{lemmas34}.
\begin{lemma}\label{lemmas34}
{Suppose that  Assumption (A1) and $|\sigma|<\min\{\frac{1}{\sqrt{\mu_{c}}},\sqrt{\frac{2(c-1)}{3}}\}$ hold.}
  Then, for any initial value $z_0\in L^2_{\mathcal{F}_0}(\Omega;L^2(0,1))$, the solution of the closed-loop system  \dref{SPDE-Heat-trans-U0-closed} satisfies
\begin{equation}
\mathbb{E}|z(\cdot,t)|_{L^2(0,1)}^2\leq \Gamma e^{{-\theta^{*}t}},\; \forall t\geq 0,
\end{equation}
where the positive constant $\Gamma$  specified in (\ref{3603}) depends on the initial value $z_{0}$ { and
the $\theta$ for the second case of (\ref{canostff}).}
\end{lemma}
{\bf Proof.} See {\bf ``Proof of Lemma  \ref{lemmas34}''} in Appendix C.

 
The resulting closed-loop system comprised of  \dref{SPDE-Heat}, \dref{dobc-part} and \dref{controllerde} is
then written as
\begin{equation}\label{SPDE-Heat-closed}
 \left\{\begin{array}{l}
dy(x,t)=y_{xx}(x,t)dt+a(x)y(x,t)dt+\sigma y(x,t)dB(t),\;\cr 
y_x(0,t)=0,\; t\geq0,\cr y_x(1,t)=k(1,1)y(1,t)+\int_0^1k_x(1,\zeta)y(\zeta,t)d\zeta\cr\hspace{1.3cm}+w(t)-\widehat{w}(t),\; t\geq0,\cr
z(x,t)=y(x,t)-\int_0^xk(x,\zeta)y(\zeta,t)d\zeta,\; t\geq 0, \cr Z(t)=\int_0^1\cos(\pi x) z(x,t)dx,\; t\geq0,\cr
\dot{\vartheta}(t)=A\vartheta(t)+[AL+(c+\pi^{2})L]Z(t),\; t>0,\cr
\widehat{\xi}(t)=\vartheta(t)+LZ(t),\ \ \widehat{w}(t)=C\widehat{\xi}(t),\; t\geq 0.
\end{array}\right.
\end{equation}

The mean square exponential stability of the closed-loop system (\ref{SPDE-Heat-closed}) is summarized in the following
Theorem \ref{Thm-closed-spde}.

\begin{theorem} \label{Thm-closed-spde}
  {Suppose that  Assumption (A1) and $|\sigma|<\min\{\frac{1}{\sqrt{\mu_{c}}},\sqrt{\frac{2(c-1)}{3}}\}$ hold.}
 Then,  for any initial value $y_0\in L^2_{\mathcal{F}_0}(\Omega;L^2(0,1))$,
 the closed-loop system \dref{SPDE-Heat-closed} admits a unique {weak} solution
$y\in C_{\mathbb{F}}(0,+\infty;L^2(\Omega; L^2(0,1)))$; Moreover, the solution of the closed-loop system \dref{SPDE-Heat-closed}
satisfies
\begin{equation}
\mathbb{E}|y(\cdot,t)|^{2}_{L^2(0,1)}\leq \Gamma^{*}e^{{-\theta^{*}t}},\; \forall t\geq 0,
\end{equation}
where $\Gamma^{*}$ specified in (\ref{constatsd}) is a positive constant depending  on the initial value $y_0$.
\end{theorem}

{\bf Proof.}  See {\bf ``Proof of Theorem \ref{Thm-closed-spde}''} in  Appendix D.


The almost surely exponential stability of the closed-loop system (\ref{SPDE-Heat-closed}) is summarized in the following
Theorem \ref{fde67g}.

\begin{theorem}\label{fde67g}
{Suppose that  Assumption (A1)  and $|\sigma|<\min\{\frac{1}{\sqrt{\mu_{c}}},\sqrt{\frac{2(c-1)}{3}}\}$ hold.}
 Then, for any initial value $y_0\in L^2_{\mathcal{F}_0}(\Omega;L^2(0,1))$,
 the closed-loop system \dref{SPDE-Heat-closed} admits a unique {weak} solution
$y\in C_{\mathbb{F}}(0,+\infty;L^2(\Omega; L^2(0,1)))$; Moreover, the solution of the closed-loop system \dref{SPDE-Heat-closed} satisfies
\begin{equation}
\limsup_{t\to \infty}\frac{1}{t}\log|y(\cdot,t)|_{L^2(0,1)} \leq 
{-\frac{\theta^{*}}{2}} \;\;\;\; \mbox{almost surely.}
\end{equation}
\end{theorem}

{\bf Proof.} See {\bf ``Proof of Theorem \ref{fde67g}''} in  Appendix E.

{\begin{remark} It is noted that
 system (1) contains two classes of disturbances: one is the in-domain multiplicative noise $\sigma y(x,t)dB(t)$
and the other is the unknown boundary external disturbance $w(t)$. The in-domain multiplicative noise whose
 intensity $\sigma$ satisfying $|\sigma|<\min\{\frac{1}{\sqrt{\mu_{c}}},\sqrt{\frac{2(c-1)}{3}}\}$
 is attenuated by the proposed DOBC boundary control \dref{controllerde} in a passive way. However,  $w(t)$
is rejected by  an active  estimation/cancellation strategy in the feedback loop.  
  \end{remark}}

{\begin{remark}\label{rek-Lc5}
For practical examples with the known matrices $A,L,C$ and the chosen parameter $c$,
it is  seen from (\ref{ucontant}) that
 $\mu_{c}$ can be computed directly so that the inequality condition $|\sigma|<\min\{\frac{1}{\sqrt{\mu_{c}}},\sqrt{\frac{2(c-1)}{3}}\}$
 in Lemma \ref{lemmas34}, Theorems \ref{Thm-closed-spde} and \ref{fde67g} can be checked; And for given known matrices $A,L,C$,
 the relationship of $u_{c}$ with the chosen parameter $c$ can also be
determined, which will provide a guideline on how to determine the parameter $c$ to satisfy the inequality condition.
Significantly, we can also see that the
maximum tolerance of the in-domain multiplicative noise is $\min\{\frac{1}{\sqrt{\mu_{c}}},\sqrt{\frac{2(c-1)}{3}}\}$
which depends on $c$ and $L$. This is reasonable because they are used in the design of both the gain $(c+\pi^{2})L$ of the disturbance observer (\ref{dobc-part})
and the function gain  $k(\cdot,\cdot)$ of the DOBC boundary control (\ref{controllerde}).
\end{remark}}
{\begin{remark}
From Remark \ref{rek-Lc5}, we know that the maximum tolerance of the multiplicative noise depends on the designed $c$ and $L$.
For practical example, it is possible to guarantee that the maximum tolerance can be sufficiently large by appropriately designing $c$ and $L$.
For example, when  $n=2$, $A=\sbm{0& 2\\ -2& 0}$ and $C=[1,0]$, the  matrix $L=[l_1,  l_2]^\top$ where $l_1<0$ and $l_2<2$ is chosen so that
$A+LC=\sbm{l_1& 2\\ l_2-2& 0}$ is Hurwitz. Moreover,
from \dref{Ly-alg-equ}-\dref{ucontant}, it is seen that $\mu_{c}$ is calculated to be $\mu_{c}=\lambda_{\max}\left(\sbm{1 &0\\  L&0 }^{\top}Q_{c}\sbm{1&0\\  L&0 }\right)=\frac{1}{2c+2\pi^2}+\frac{-l_1(c+\pi^2)-l_2}{(c+\pi^2)\lambda^{*}}+(l_2-3+\frac{2c+2\pi^2}{\lambda^{*}})\frac{l_1}{2}-\frac{l_1l_2}{2}
+\frac{1}{l_2-2}[\frac{l_2^2}{\lambda^{*}}+\frac{l_1l_2^2}{2}- (l_2-3+\frac{2c+2\pi^2}{\lambda^{*}})\frac{l_2^2}{l_1}]>0$, where $\lambda^{*}=(l_1-c-\pi^2)(c+\pi^2)+l_2-1$.
Letting $l_2=l_1$ with $|l_1|<2$, it is clear that $\mu_c\to0$ as $(c,l_2)\to(\infty,0^-)$. Hence,
$\min\{\frac{1}{\sqrt{\mu_{c}}},\sqrt{\frac{2(c-1)}{3}}\}\to\infty$ as $(c,l_2)\to(\infty,0^-)$, which implies that
the maximum tolerance $\min\{\frac{1}{\sqrt{\mu_{c}}},\sqrt{\frac{2(c-1)}{3}}\}$ can be as large as possible by choosing properly both $c$ and $L$.
Certainly, finding  $c$ and $L$ to ensure $\min\{\frac{1}{\sqrt{\mu_{c}}},\sqrt{\frac{2(c-1)}{3}}\}$ to be large enough is generally
complicated when $n>2$.
\end{remark}
}

{\begin{remark}
Compared with the deterministic setting, many new concerns should be addressed for the boundary feedback stabilization
of stochastic PDEs subject to boundary external disturbance. On the one hand, the admissibility property
that the It\^{o} integral $\int_0^te^{\mathcal{A}(t-s)}\mathcal{B}u(s)dB(s)\in L^2_{\mathbb{F}}(\Omega; L^2(0,1))$ for any fixed $t$
 in stochastic setting  is difficult to be proved, where $e^{\mathcal{A}t}$ and $\mathcal{B}=\delta(x-1)$ are the involved semigroup and boundary control operator, respectively.
Thus, the proof of Lemma \ref{lem-existencewell-u-sys} is divided into three steps, and some new approaches including stochastic analysis to obtain the martingale property and the approximating method
are adopted to prove the well-posedness of the closed-loop system (\ref{SPDE-Heat-trans-U0-closed}).
On the other hand, compared with the deterministic situation,  it is not only a  simple application of
 the It\^{o}'s differentiation rule in stochastic setting.
For example, in the proof of Lemma \ref{lem-existencewell-u-sys},  the stopping time technique should be  used
because no available results to exclude directly the possibility that $|z(\cdot,t)|_{L^2(0,1)}^2\rightarrow\infty$ in finite time so that even
$\mathbb{E}|z(\cdot,t)|^2_{L^2(0,1)}$ does not necessarily exist and $\int^{t}_{0}\sigma|z(\cdot,s)|_{L^2(0,1)}^2dB(t)$ is not necessarily a martingale.
In addition, very few references address the almost surely boundary stabilization of stochastic PDEs with or without boundary external disturbance.
Thus, in the proof of Theorem \ref{fde67g},  some stochastic analysis techniques including Chebyshev's inequality,
 Burkholder-Davis-Gundy inequality and Borel-Cantelli lemma are used for the proof of almost surely exponential stability of the closed-loop system \dref{SPDE-Heat-closed}, where It\^{o}'s differentiation rule is just the first step.
\end{remark}}

\section{Numerical simulations}\label{Se4}

\begin{figure}[ht]\centering
\subfigure[The open-loop state $y(x,t)$]
 {\includegraphics[width=7.3cm,height=6.1cm]{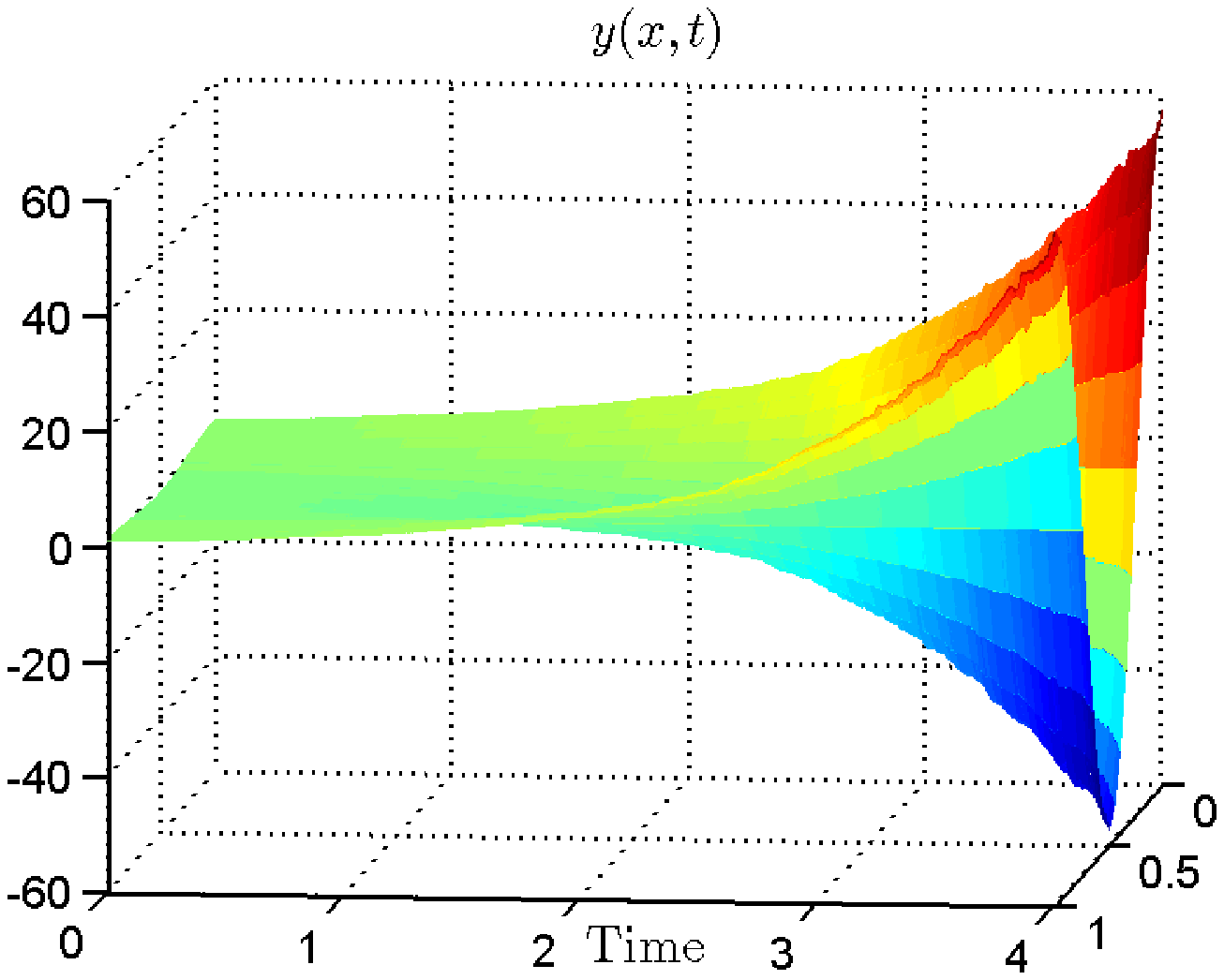}\label{Fig-ref-con-open}}
\subfigure[The closed-loop state $y(x,t)$]
 {\includegraphics[width=7.3cm,height=6.1cm]{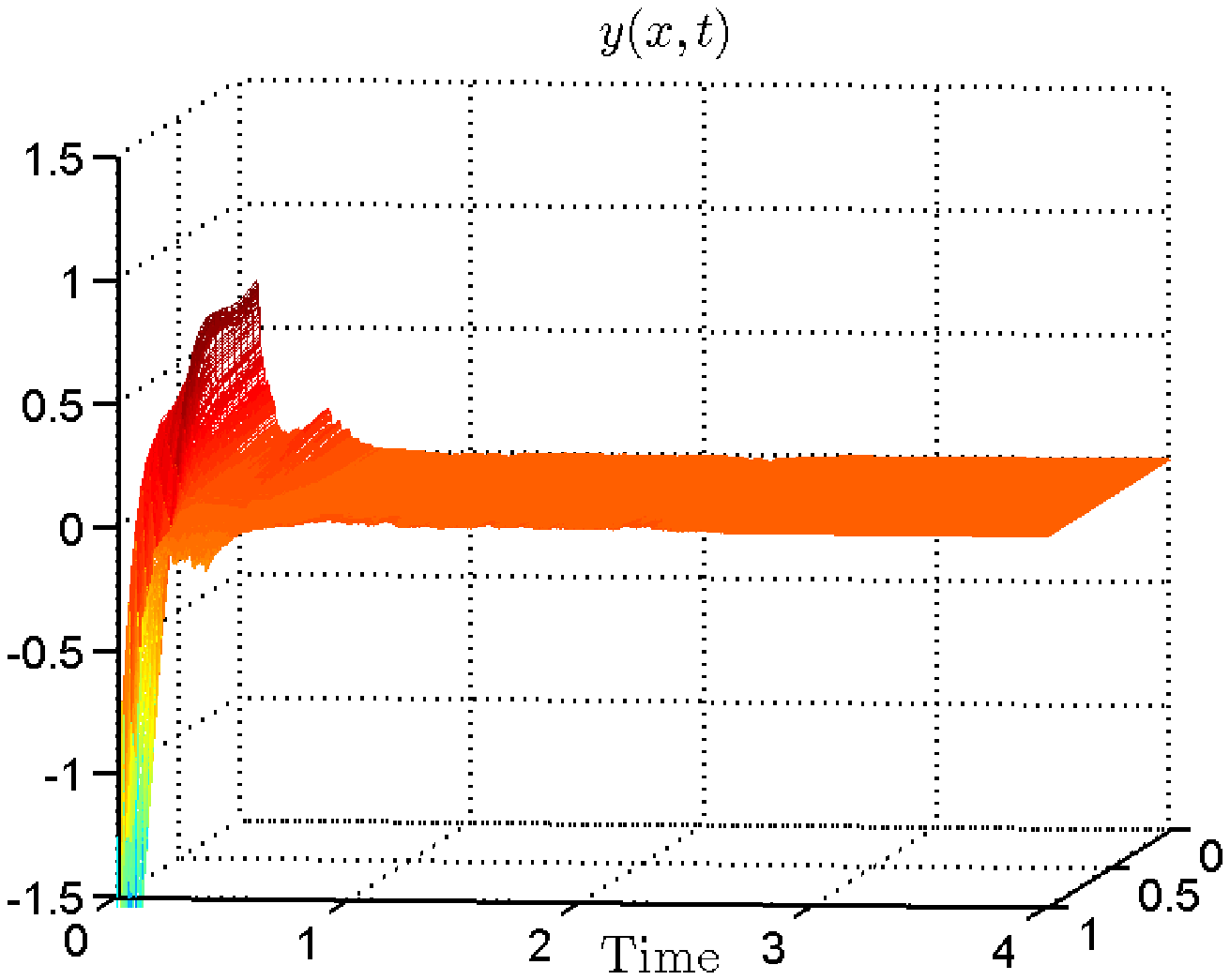}\label{Fig-con-b}}
\caption{The state $y(x,t)$ of the open-loop system without boundary control and the
state $y(x,t)$ of the closed-loop system \dref{SPDE-Heat-closed}.\;(for
interpretation of the references to color of the figure's legend in
this section, we refer to the PDF version of this paper).}
\label{Fig-ref-con}
\end{figure}

In this section, we present some numerical simulations for the closed-loop system \dref{SPDE-Heat-closed}
  for illustration of the effectiveness of the proposed DOBC control.  As pointed out in Remark \ref{Rem-exam},
  we take $a(x)=4\pi^2+1.005$, $\sigma=0.1$, $y_0(x)=\cos(2\pi x)$
in system \dref{SPDE-Heat}. In this case,  system  \dref{SPDE-Heat} is not stable if there is no boundary control input.
{This can be seen from Figure \ref{Fig-ref-con-open}.}

The known matrices in the exogenous system
  \dref{exogenous} are chosen  as $A=\sbm{0& 2\\ -2& 0}$ and  $C=[1,0]$, and
$L$ is taken as $L=[-5,-1]^\top$ to guarantee that $A+LC=\sbm{-5&2\\ -3&0}$ is Hurwitz.
The parameter $c$  in system \dref{SPDE-Heat-closed} is chosen as $c=1.02$.
Solve \dref{Ly-alg-equ} to get
$
{Q_{c}}=\begin{bmatrix}
0.0134 &   0.0041&    0.0041\\ 0.0041  &  0.1667  &  0.1667\\ 0.0041 &   0.1667  &  0.6667
\end{bmatrix}
$
and $\lambda_{\max}({Q_{c}})=0.7172$. Calculate \dref{ucontant} to get ${\mu_{c}}=6.464$.
Clearly, the conditions in Lemma \ref{lem-Z-e-to0}, Theorem \ref{Thm-closed-spde} and Theorem \ref{fde67g} are satisfied.

\begin{figure}[ht]\centering
 \subfigure[$(\widehat{w}(t),w(t),\widehat{w}(t)-w(t))$]
 {\includegraphics[width=7.3cm,height=6.1cm]{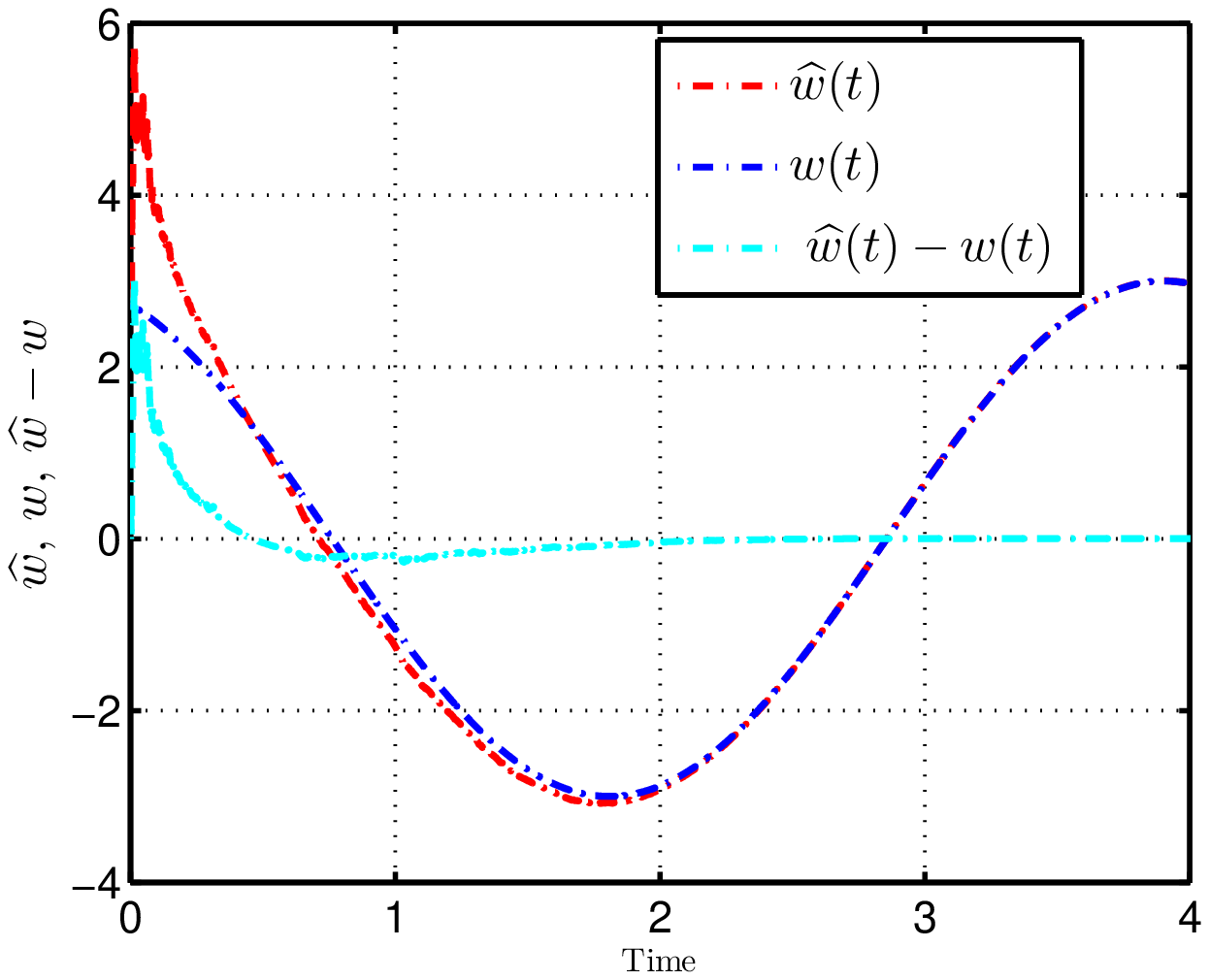}\label{Fig-refhatw-c}}
  \subfigure[ $(u(t),u(t)+\hat{w}(t))$]
 {\includegraphics[width=7.3cm,height=6.1cm]{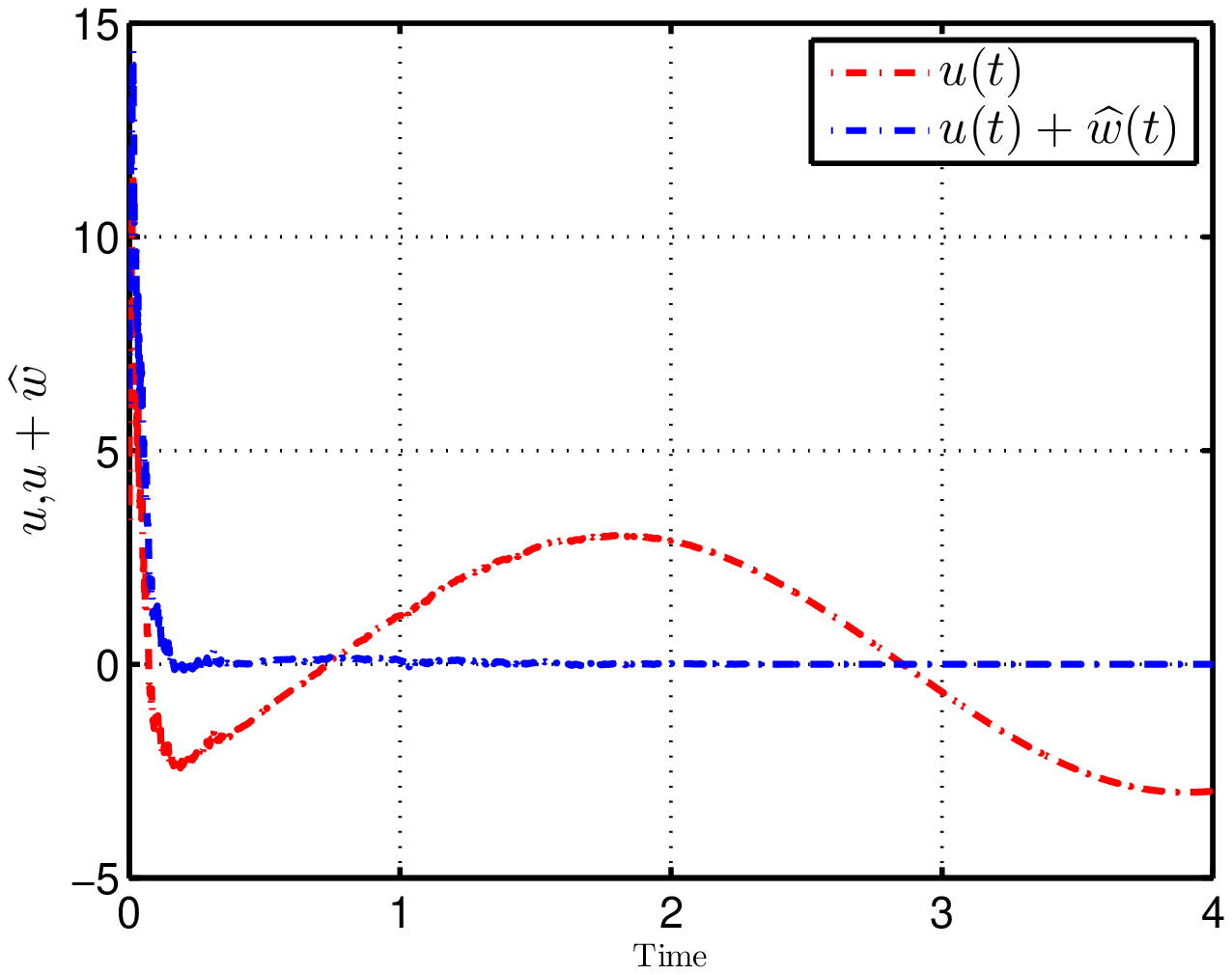}\label{Fig-refeta-ccc}}
\caption{The estimation of boundary external disturbance $(\widehat{w}(t),w(t),\widehat{w}(t)-w(t))$,
DOBC boundary control $u(t)$ and feedback control $u(t)+\hat{w}(t)$ \;(for
interpretation of the references to color of the figure's legend in
this section, we refer to the PDF version of this paper).}
\label{Fig-ref-con-vec}
\end{figure}
It is seen from Figure \ref{Fig-con-b}  that the closed-loop state $y(x,t)$ is convergent to zero quite quickly.
It is observed from Figure \ref{Fig-ref-con-vec}(a) that the boundary external disturbance $w(t)$ is effectively estimated by $\widehat{w}(t)$.
{Finally, the effect of both DOBC boundary control $u(t)$ and feedback control $u(t)+\hat{w}(t)$ can be seen from Figure \ref{Fig-refeta-ccc}.}

\section{Concluding remarks}\label{Se5}
This paper is the first effort  to apply the disturbance observer-based control (DOBC) approach to
the stabilization and disturbance rejection for the  stochastic distributed parameter systems.
The considered plant is a one-dimensional anti-stable stochastic heat equation
driven by multiplicative white noise subject to boundary external  disturbance
generated from an exogenous system. A boundary feedback stabilizing control
is first designed for the system without external boundary
disturbance by the backstepping approach. A disturbance observer is designed to estimate in real time
the boundary external  disturbance, and then a DOBC boundary control
constructed by the boundary feedback stabilizing control and
 a compensation term by use of the estimate of the
boundary external  disturbance is designed. It is shown  that
the resulting closed-loop system is exponentially stable in mean square and
almost surely. Some numerical simulations validate the
theoretic results and effectiveness of the proposed DOBC approach.
{A further development of this topic is to address the output feedback boundary
 control for stochastic PDEs subject to unknown external disturbances
only by pointwise measurement.}

\vspace{0.2cm}

{\bf APPENDIX A: Proof of Lemma \ref{lem-Z-e-to0}.}
Let $V(Z,\eta)=\sbm{Z\\ \eta}^{\top}{Q_{c}}\sbm{Z\\ \eta}$. {By \dref{coupled-Z-e-sys} and \dref{Ly-alg-equ}, and applying It\^{o}'s formula to $V(Z(t),\eta(t))$ with respect to $t$,
we obtain}
{\begin{eqnarray}
\hspace{-0.6cm}&& \mathbb{E}V(Z(t),\eta(t))=\mathbb{E}V(Z(0),\eta(0))-\int^{t}_{0}\mathbb{E}[|Z(s)|^2+|\eta(s)|^2]ds\cr \hspace{-0.6cm}&&
+\sigma^{2}\mathbb{E}\int^{t}_{0}\sbm{Z(s)\\ \eta(s)}^{\top}\sbm{1 &0\\  L&0 }^{\top}{Q_{c}}\sbm{1&0\\  L&0 }\sbm{Z(s)\\ \eta(s)}ds.
\end{eqnarray}}
It is also easily obtained that
\begin{equation}\label{equatilytdd}
\lambda_{\min}({Q_{c}})[|Z|^2+|\eta|^2]\leq V(Z,\eta)\leq \lambda_{\max}({Q_{c}})[|Z|^2+|\eta|^2].
\end{equation}
In this way,
\begin{eqnarray}\label{3277}
 \frac{d}{dt}\mathbb{E}V(Z(t),\eta(t))\leq
-\frac{(1-{\mu_{c}}\sigma^2)}{\lambda_{\max}({Q_{c}})}\mathbb{E}V(Z(t),\eta(t)),
\end{eqnarray}
where $\mu_{c}$ is defined as that in \dref{ucontant}. The \dref{estiomerror} can follow from (\ref{equatilytdd}) and (\ref{3277}) that
\begin{eqnarray}
&&\mathbb{E}[|Z(t)|^2+|\eta(t)|^2] \cr
&&\leq\frac{\mathbb{E}V(Z(t),\eta(t))}{\lambda_{\min}({Q_{c}})}\leq
\frac{\mathbb{E}V(Z(0),\eta(0))}{\lambda_{\min}({Q_{c}})}e^{-\frac{(1-{\mu_{c}}\sigma^2)}{\lambda_{\max}({Q_{c}})}t}   \cr &&
\leq\frac{ \lambda_{\max}({Q_{c}})}{\lambda_{\min}({Q_{c}})}\mathbb{E}[|Z(0)|^2+|\eta(0)|^2]e^{-\frac{(1-{\mu_{c}}\sigma^2)}{\lambda_{\max}({Q_{c}})}t}.
\end{eqnarray}
Since both the drift term and diffusion term of system (\ref{coupled-Z-e-sys}) satisfy the linear growth condition, by
\cite[Theorem 4.2, p. 128]{mao2007}, \dref{almostlysure}
can be directly concluded from \dref{estiomerror}, where the relevant parameters are specified as  $\lambda=\frac{(1-{\mu_{c}}\sigma^2)}{\lambda_{\max}({Q_{c}})}$, $p=2$.
This completes the proof of the lemma. \hfill$\Box$

\vspace{0.2cm}

{\bf APPENDIX B: Proof of Lemma \ref{lem-existencewell-u-sys}.}
{ Since the admissibility theory  for stochastic setting is not available for our case, the well-posedness cannot be simply proved.
To derive the well-posedness of \dref{SPDE-Heat-trans-U0-closed}, we split the proof into  three
steps.}

{\bf Step 1:  Suppose that $\widetilde{w}\in
C^1_{\mathbb{F}}(0,T;L^2(\Omega;\mathbb{R}))$.}  Consider the
boundary value problem of a second order PDE
subject to random boundary value $\widetilde{w}(t)$ as follows:
\begin{equation}
\left\{
\begin{array}{l}
p_{xx}(x,t)-cp(x,t)=0,\ x\in(0,1),\cr
p_x(0,t)=0,\ p_x(1,t)=\widetilde{w}(t).
\end{array}
\right.
\end{equation}

We first claim that $p\in C^1_{\mathbb{F}}(0,T;L^2(\Omega;H^1(0,1)))$.
Actually, by the classical theory of elliptic equations with Neumann
boundary condition, for any $T>0$, there exists a positive constant $C(T)$ depending  on $T$ only, such that
\begin{equation}
|p|_{H^1(0,1)}^2\leq C(T)|\widetilde{w}|^2,\; \forall\;
t\in[0,T] \;\;\;\mbox{almost surely.}
\end{equation}
This, together with the fact   $\widetilde{w}(t)\in
L^2_{\mathcal{F}_t}(\Omega;\mathbb{R})$ for all $t\geq 0$ which is  concluded from (\ref{estiomerror}) and $\widetilde{w}(t)=-C\eta(t)$, implies that $p(\cdot,t)\in
L^2_{\mathcal{F}_t}(\Omega;H^1(0,1))$, and
\begin{equation}\label{v-continous}
\mathbb{E}|p(\cdot,t)-p(\cdot,s)|^2_{H^1(0,1)}\leq
C(T)\mathbb{E}|\widetilde{w}(t)-\widetilde{w}(s)|^2.
\end{equation}
Since $\widetilde{w}\in  C^1_{\mathbb{F}}(0,T;L^2(\Omega;\mathbb{R}))$, it
follows  that $p\in C^1_{\mathbb{F}}(0,T;L^2(\Omega;H^1(0,1)))$.
Next, we consider the following stochastic heat equation:
\begin{equation}\label{SPDE-Heat-trans-closed-err-v}
\left\{\begin{array}{l}
d{v}(x,t)=[{v}_{xx}(x,t)-c{v}(x,t)-p_t(x,t)]dt \cr\hspace{1.5cm}
                 +{\sigma}(v(x,t)+ p(x,t))dB(t), \cr
{v}_x(0,t)=0, \ \ \ 
{v}_x(1,t)=0,\ t\geq0,\cr
{v}(x,0)=z_0(x)-p(x,0),\;0\leq x\leq1,
\end{array}\right.
\end{equation}
which  can be rewritten as
\begin{equation}\label{absAsmc}
d{v}(\cdot,t)=[\mathcal{A}{v}(\cdot,t)-p_t(\cdot,t)]dt+\sigma[{v}(\cdot,t)+p(\cdot,t)]dB(t),
\end{equation}
where the linear operator $\mathcal{A}$ is defined  by
\begin{equation}
\left\{\begin{array}{l} [\mathcal{A}f](x)=f''(x)-cf(x), \cr
D(\mathcal{A})=\{f\in H^2(0,1)|\; f'(0)=f'(1)=0\}.
\end{array}\right.
\end{equation}

It is well known that $\mathcal{A}$ generates an exponentially  stable
$C_0$-semigroup on $L^2(0,1)$ with the decay rate $-c$. 
This, together with
\cite[Theorem 3.5]{SPDEbook2011}, shows that  system
\dref{SPDE-Heat-trans-closed-err-v} admits a unique weak solution
${v}\in C_{\mathbb{F}}(0,+\infty;L^2(\Omega; L^2(0,1)))$. Define
$z(x,t)={v}(x,t)+p(x,t)$. It is easy to see that
$z(x,t)$ is a unique weak solution to \dref{SPDE-Heat-trans-U0-closed}.

{\bf Step 2:  Suppose $\widetilde{w}\in
C^1_{\mathbb{F}}(0,T;L^2(\Omega;\mathbb{R}))$ and establish an energy estimate for the solution of
\dref{SPDE-Heat-trans-U0-closed}.}  By virtue of It\^{o}'s formula,  $dz^2(x,t)= (\sigma^2-2c)z^2(x,t)dt+2z_{xx}(x,t)z(x,t)dt+2\sigma z^2(x,t)dB(t)$
and hence
\begin{equation}\label{u-ito-dif}
\begin{array}{ll}
&d(|z(\cdot,t)|^2_{L^2(0,1)})=
(\sigma^2-2c)|z(\cdot,t)|^2_{L^2(0,1)}dt-
2|z_x(\cdot,t)|^2_{L^2(0,1)}dt\cr&\ \ \ +2z(1,t)\widetilde{w}(t)dt+2\sigma|z(\cdot,t)|_{L^2(0,1)}^2dB(t).
\end{array}
\end{equation}

For every integer $n\geq 1$, define the stopping time
{\begin{equation}\label{340h}
\tau_{n}=T\wedge\inf\{t\geq 0:|z(\cdot,t)|_{L^2(0,1)}^2\geq n\},
\end{equation}
and set $\inf\emptyset =\infty$.
Clearly, $\tau_{n}\uparrow \tau_{\infty}$ almost surely for some random variable $\tau_{\infty}$. We will show
that $\tau_{\infty}=T$ almost surely.}
By the definition of the stopping time in (\ref{340h}), $\int^{t\wedge \tau_{n}}_{0}2\sigma|z(\cdot,s)|_{L^2(0,1)}^2dB(s)$ is a martingale
for all  $t\in [0,T]$.
 We then have
\begin{equation}\label{GWl1}
\begin{array}{l}
\mathbb{E}|z(\cdot,t\wedge \tau_{n})|^2_{L^2(0,1)}=\mathbb{E}|z(\cdot,0)|^2_{L^2(0,1)}+2\mathbb{E}\int_0^{t\wedge \tau_{n}}z(1,s)\widetilde{w}(s)ds
\cr -2\mathbb{E}\int_0^{t\wedge \tau_{n}}|z_x(\cdot,s)|^2_{L^2(0,1)}ds
+\mathbb{E}\int_0^{t\wedge \tau_{n}}(\sigma^2-2c)|z(\cdot,s)|^2_{L^2(0,1)}ds.
\end{array}
\end{equation}
By the Cauchy-Schwarz inequality, Young's inequality and the
classical embedding theorem in Sobolev spaces, we obtain
\begin{eqnarray}\label{GWl2}
&&\mathbb{E}\int_0^{t\wedge \tau_{n}}z(1,s)\widetilde{w}(s)ds
 \leq \frac{1}{2\varepsilon}\mathbb{E}\int_0^{t\wedge \tau_{n}}\widetilde{w}^2(s)ds
\cr &&+\varepsilon\mathbb{E}\int_0^{t\wedge \tau_{n}}[|z(\cdot,s)|_{L^2(0,1)}^2+|z_x(\cdot,s)|_{L^2(0,1)}^2]ds,
\end{eqnarray}
where $0<\varepsilon<1$. The following proof is divided into two cases.
\begin{itemize}

\item {\it Case 1:  $\sigma^2-2c+2\varepsilon \leq 0$.}  In this case,  by
\dref{GWl1} and \dref{GWl2}, it follows that
\begin{eqnarray}\label{fd3ff}
\mathbb{E}|z(\cdot,t\wedge\tau_{n})|^2_{L^2(0,1)}
\leq\mathbb{E}|z(\cdot,0)|^2_{L^2(0,1)}+\frac{1}{\varepsilon}\mathbb{E}\int_0^{t\wedge \tau_{n}}\widetilde{w}^2(s)ds.
\end{eqnarray}
{Thus, $nP\{\tau_{n}\leq t\}\leq\mathbb{E}|z(\cdot,0)|^2_{L^2(0,1)}+\frac{1}{\varepsilon}\mathbb{E}\int_0^{T}\widetilde{w}^2(s)ds$ for any $t\in [0,T]$.
Letting $n\rightarrow\infty$, we have $P\{\tau_{\infty}\leq t\}=0\Leftrightarrow P\{\tau_{\infty }> t\}=1$ for any $t\in [0,T]$. It can be concluded
that $\tau_{\infty}=T$ almost surely and then $\tau_{n}\uparrow T$ almost surely.}
By Fatou's Lemma and passing to the limit as  $n\rightarrow\infty$ {for \dref{fd3ff}}, we have,  for all $t\in [0,T]$, that
\begin{eqnarray}
\mathbb{E}|z(\cdot,t)|^2_{L^2(0,1)}\leq
\mathbb{E}|z(\cdot,0)|^2_{L^2(0,1)}+\frac{1}{\varepsilon}\int_0^T\mathbb{E}\widetilde{w}^2(s)ds,
\end{eqnarray}
which implies \dref{u-theta-estinq};

\item {\it Case 2:  $\sigma^2-2c+2\varepsilon > 0$.} In this case, {similarly to case 1, $\tau_{n}\uparrow T$ almost surely,}  and
by Fatou's Lemma and passing to the limit as  $n\rightarrow\infty$, we can also obtain
 for all $t\in [0,T]$ that
\begin{equation}
\begin{array}{l}
 \mathbb{E}|z(\cdot,t)|^2_{L^2(0,1)}  \leq
\mathbb{E}|z(\cdot,0)|^2_{L^2(0,1)} +\frac{1}{\varepsilon}\int_0^{T}\mathbb{E}\widetilde{w}^2(s)ds
\cr
+(\sigma^2-2c+2\varepsilon)\int_0^{t}\mathbb{E}|z(\cdot,s)|^2_{L^2(0,1)}ds.
\end{array}
\end{equation}
By Gronwall's inequality, for all $t\in [0,T]$, we have
\begin{eqnarray}
&&\hspace{-1cm}\mathbb{E}|z(\cdot,t)|^2_{L^2(0,1)}\cr&&\hspace{-1cm}\leq (\mathbb{E}|z(\cdot,0)|^2_{L^2(0,1)}+\frac{1}{\varepsilon}\int_0^{T}\mathbb{E}\widetilde{w}^2(s)ds) e^{(\sigma^2-2c+2\varepsilon)t},
\end{eqnarray}
which implies \dref{u-theta-estinq}.
\end{itemize}

{\bf Step 3:} It follows from (\ref{estiomerror}) in  Lemma \ref{lem-Z-e-to0}
that $\widetilde{w}\in L^2_{\mathbb{F}}(0,T;L^2(\Omega;\mathbb{R}))$.  We  can then find a sequence
$\{\widetilde{w}_n\}_{n=1}^\infty\subset
C^1_{\mathbb{F}}(0,T;L^2(\Omega;\mathbb{R}))$ such that $\lim_{n\to\infty}\widetilde{w}_n=\widetilde{w} \mbox{ in } L^2_{\mathbb{F}}(0,T;L^2(\Omega;\mathbb{R})).$
Denote by $z_n(x,t)$ the {weak} solution to \dref{SPDE-Heat-trans-U0-closed}
with the initial value  $z_0(x)$ and  random boundary value
$\widetilde{w}_n(t)$. Then,   $\{z_n(x,t)\}_{n=1}^\infty$ is a
Cauchy sequence in $C_{\mathbb{F}}([0,T];L^2(\Omega;L^2(0,1)))\cap
L^2_{\mathbb{F}}(0,T;L^2(\Omega;H^1(0,1)))$. Thus, there exist a
unique $z\in C_{\mathbb{F}}([0,T];L^2(\Omega;L^2(0,1)))\cap
L^2_{\mathbb{F}}(0,T;L^2(\Omega;H^1(0,1)))$ such that $\lim_{n\to\infty}z_n=z\mbox{ in } C_{\mathbb{F}}([0,T];L^2(\Omega;L^2(0,1)))\cap L^2_{\mathbb{F}}(0,T; $ $L^2(\Omega;H^1(0,1))).$
From the definition of $z_n(x,t)$, we have,   for all $t\in[0,T]$
and $\phi\in H^1(0,1)$, that
\begin{eqnarray}\label{un-equ}
\hspace{-1.2cm}&&\int_0^1z_n(x,t)\phi(x)dx-\int_0^1z_n(x,0)\phi(x)dx\cr \hspace{-1.2cm} &&=\int_0^t\frac{\partial z_{n}(1,s)}{\partial x}\phi(1)ds
-\int_0^t\int_0^1 \frac{z_{n}(x,s)}{\partial x}\phi'(x)dxds\cr \hspace{-1.2cm} &&-c\int_0^t\int_0^1z_n(x,s)\phi(x)dxds+\int_0^t\int_0^1\sigma
z_n(x,s)\phi(x)dxdB(s).
\end{eqnarray}
This yields, for all $t\in[0,T]$, that
\begin{eqnarray}\label{u-equ}
\hspace{-1cm}&& \int_0^1z(x,t)\phi(x)dx-\int_0^1z(x,0)\phi(x)dx \cr\hspace{-1cm} &&=\int_0^tz_x(1,s)\phi(1)ds
-\int_0^t\int_0^1z_x(x,s)\phi'(x)dxds\cr\hspace{-1cm} &&-c\int_0^t\int_0^1z(x,s)\phi(x)dxds+\int_0^t\int_0^1\sigma
z(x,s)\phi(x)dxdB(s).
\end{eqnarray}
Therefore,  $z(x,t)$ is a {weak} solution to
\dref{SPDE-Heat-trans-U0-closed} and satisfies \dref{u-theta-estinq}. \hfill$\Box$

\vspace{0.2cm}

{\bf  APPENDIX C: Proof of Lemma \ref{lemmas34}.} Let $\beta(x,t)=z(x,t)+\frac{x^{2}}{2}C\eta(t).$
Clearly,
$d\beta(x,t)=dz(x,t)+\frac{x^2}{2}Cd\eta(t)$
and $\beta_{xx}(x,t)=z_{xx}(x,t)+C\eta(t)$.
By \dref{e-sde-F}, a direct computation shows that $\beta(x,t)$ satisfies the following stochastic PDE:
\begin{equation}\label{beta-spde-II}
 \left\{
\begin{array}{l}
d\beta(x,t)\!\!=[\beta_{xx}(x,t)-c\beta(x,t)+(\frac{cx^2}{2}-1)C\eta(t)
+\frac{x^2}{2}C(A+LC)\eta(t)]dt\crr\hspace{1.4cm}
+\sigma[\frac{x^{2}}{2}CLZ(t)+\beta(x,t)-\frac{x^2}{2}C\eta(t)]dB(t), \cr
\beta_x(0,t)=\beta_x(1,t)=0,\ \ t\geq0.
\end{array}
\right.
\end{equation}

By It\^{o}'s formula, a direct computation shows that
\begin{eqnarray}\label{3522}
\hspace{-0.6cm}&&  d\beta^2(x,t)\cr\hspace{-0.6cm}  &&=2\beta(x,t)d\beta(x,t)+\sigma^{2}[\frac{x^{2}}{2}CLZ(t)+\beta(x,t)-\frac{x^2}{2}C\eta(t)]^2dt\cr \hspace{-0.6cm} &&
=2\beta(x,t)[\beta_{xx}(x,t)-c\beta(x,t)+(\frac{cx^2}{2}-1)C\eta(t)+\frac{x^2}{2}C\cdot\cr\hspace{-0.6cm}  &&(A+LC)\eta(t)]dt
+\sigma^{2}[\frac{x^{2}}{2}CLZ(t)+\beta(x,t)-\frac{x^2}{2}C\eta(t)]^2dt\cr\hspace{-0.6cm}  &&
+2\sigma\beta(x,t)[\frac{x^2}{2}C LZ(t)+\beta(x,t)-\frac{x^2}{2}C\eta(t)]dB(t).
\end{eqnarray}
 Since from Lemmas \ref{lem-Z-e-to0} and  \ref{lem-existencewell-u-sys},  $2\sigma\int_0^t\int_0^1\beta(x,s)[\frac{x^2}{2}C LZ(s)+\beta(x,s)-\frac{x^2}{2}C\eta(s)]dxdB(s)$
 is a martingale for all $t\geq 0$. Integrating on both sides of (\ref{3522}) with respect to   $x$ and $t$ and
 taking mathematical expectations, we obtain
\begin{eqnarray}
 \hspace{-0.8cm}&&\mathbb{E}\int_0^1\beta^2(x,t)dx =\mathbb{E}\int_0^1\beta^2(x,0)dx-2c\int_0^t\mathbb{E}\int_0^1\beta^2(x,s)dxds\cr\hspace{-0.8cm}&&
-2\int_0^t\mathbb{E}\int_0^1\beta_{x}^2(x,s)dxds
+2\int_0^t\mathbb{E}\int_0^1\beta(x,s)\cdot\cr\hspace{-0.8cm}&&[(\frac{cx^2}{2}-1)C\eta(s)+\frac{x^2}{2}C(A+LC)\eta(s)]dxds\cr\hspace{-0.8cm}&&
+\int_0^t\mathbb{E}\int_0^1\sigma^{2}[\frac{x^{2}}{2}CLZ(s)+\beta(x,s)-\frac{x^2}{2}C\eta(s)]^2dxds,
\end{eqnarray}
which implies that
\begin{eqnarray}\label{3566}
 \hspace{-0.7cm}&&
\frac{d}{dt}\mathbb{E}\int_0^1\beta^2(x,t)dx\cr
\hspace{-0.7cm}&&=-2c\mathbb{E}\int_0^1\beta^2(x,t)dx
-2\mathbb{E}\int_0^1\beta_{x}^2(x,t)dx\cr\hspace{-0.7cm}&&
+2\mathbb{E}\int_0^1\beta(x,t)[(\frac{cx^2}{2}-1)C\eta(t)+\frac{x^2}{2}C(A+LC)\eta(t)]dx\cr\hspace{-0.7cm}&&
+\mathbb{E}\int_0^1\sigma^{2}[\frac{x^{2}}{2}CLZ(t)+\beta(x,t)-\frac{x^2}{2}C\eta(t)]^2dx\cr \hspace{-0.7cm}&&
\leq-2c\mathbb{E}\int_0^1\beta^2(x,t)dx+\mathbb{E}\int_0^1\beta^2(x,t)dx
+|(1+\frac{c}{2})C|^{2}\mathbb{E}|\eta(t)|^2\cr\hspace{-0.7cm}&&+
\mathbb{E}\int_0^1\beta^2(x,t)dx+|\frac{C(A+LC)}{2}|^{2}\mathbb{E}|\eta(t)|^2
+\frac{3\sigma^{2}|CL|^{2}}{4}\mathbb{E}|Z(t)|^{2}\cr\hspace{-0.7cm}&&+3\sigma^{2}\mathbb{E}\int_0^1\beta^2(x,t)dx+\frac{3\sigma^{2}|C|^{2}}{4}\mathbb{E}|\eta(t)|^2.
\end{eqnarray}
By (\ref{estiomerror}) in Lemma \ref{lem-Z-e-to0} and (\ref{3566}),  we conclude that
{\begin{eqnarray}\label{359fd}
 \hspace{-0.6cm}&&\mathbb{E}\int_0^1\beta^2(x,t)dx\cr\hspace{-0.6cm}&& 
 \leq e^{-(2c-2-3\sigma^{2})t}\mathbb{E}\int_0^1\beta^2(x,0)dx\cr\hspace{-0.6cm}&&
+\Gamma_{1}\int_0^te^{-(2c-2-3\sigma^{2})(t-s)}\mathbb{E}[|\eta(s)|^2+|Z(s)|^2]ds \cr\hspace{-0.6cm}&&
\leq e^{-(2c-2-3\sigma^{2})t}\mathbb{E}\int_0^1\beta^2(x,0)dx
+\Gamma_{1} \frac{\lambda_{\max}(Q_{c})}{\lambda_{\min}(Q_{c})}\mathbb{E}[|Z(0)|^2\cr\hspace{-0.6cm}&&+|\eta(0)|^2]\int_0^te^{-(2c-2-3\sigma^{2})(t-s)}e^{-\frac{1-\mu_{c}\sigma^{2}}{\lambda_{\max}(Q_{c})}s}ds
 \leq \Gamma_{2}e^{-\theta^{*}t},
\end{eqnarray}
where $\theta^{*}$ is given in (\ref{canostff}),   $\Gamma_{1}= \max\{|(1+\frac{c}{2})C|^{2}+|\frac{C(A+LC)}{2}|^{2}+\frac{3\sigma^{2}|C|^{2}}{4},\frac{3\sigma^{2}|CL|^{2}}{4}\}$,
$\Gamma_{2}=\mathbb{E}\int_0^1\beta^2(x,0)dx+\frac{\Gamma_{1} }{|2c-2-3\sigma^{2}-\frac{1-\mu_{c}\sigma^{2}}{\lambda_{\max}(Q_{c})}|}\frac{\lambda_{\max}(Q_{c})}{\lambda_{\min}(Q_{c})}\mathbb{E}[|Z(0)|^2+|\eta(0)|^2]$
if $2c-2-3\sigma^{2}\neq\frac{1-\mu_{c}\sigma^{2}}{\lambda_{\max}(Q_{c})}$, and $\Gamma_{2}=\mathbb{E}\int_0^1\beta^2(x,0)dx+\sup_{t\geq0}te^{-(2c-2-3\sigma^{2}-\theta)t}\frac{\Gamma_{1} \lambda_{\max}(Q_{c})}{\lambda_{\min}(Q_{c})}\mathbb{E}[|Z(0)|^2+|\eta(0)|^2]$ otherwise with $\theta$ given in (\ref{canostff}).
Furthermore, from  (\ref{estiomerror}) in Lemma \ref{lem-Z-e-to0} and (\ref{359fd}), there holds
\begin{eqnarray}
\hspace{-0.7cm}&& \mathbb{E}\int_0^1 z^2(x,t)dx
\leq 2\mathbb{E}\int_0^1\beta^2(x,t)dx+\frac{|C|^{2}}{2}\mathbb{E}|\eta(t)|^{2}\cr\hspace{-0.7cm} &&
\leq2\Gamma_{2}e^{-\theta^{*}t}+\frac{|C|^{2} \lambda_{\max}(Q_{c})}{2\lambda_{\min}(Q_{c})}\mathbb{E}[|Z(0)|^2+|\eta(0)|^2]e^{-\frac{(1-\mu_{c}\sigma^2)}{\lambda_{\max}(Q_{c})}t}
\cr \hspace{-0.7cm}&& \leq\Gamma e^{-\theta^{*}t},
\end{eqnarray}
where we set
\begin{eqnarray}\label{3603}
\Gamma= 2\Gamma_{2}+\frac{|C|^{2}}{2}\frac{ \lambda_{\max}(Q_{c})}{\lambda_{\min}(Q_{c})}\mathbb{E}[|Z(0)|^2+|\eta(0)|^2].
\end{eqnarray}}
 \hfill$\Box$

{\bf APPENDIX D: Proof of Theorem \ref{Thm-closed-spde}.} The existence of the solution $y\in C_{\mathbb{F}}(0,+\infty;$ $L^2(\Omega; L^2(0,1)))$ can be concluded directly from Lemma \ref{lem-existencewell-u-sys}. In addition, it follows from (\ref{bstep-inv}), Lemma  \ref{lemmas34}  and the H\"{o}der inequality that
\begin{eqnarray}
\hspace{-0.5cm}&& \mathbb{E}\int_0^1 y^2(x,t)dx\cr\hspace{-0.5cm}&& \leq
2\mathbb{E}\int_0^1 z^2(x,t)dx+2\mathbb{E}\int_0^1\Big(\int_0^xl(x,\zeta)z(\zeta,t)d\zeta\Big)^{2}dx \cr\hspace{-0.5cm}&&
\leq 2\mathbb{E}\int_0^1 \hspace{-0.1cm}z^2(x,t)dx+2\mathbb{E}\int_0^1\hspace{-0.1cm}\Big(\hspace{-0.1cm}\int_{0}^{x}\hspace{-0.1cm}l^{2}(x,\zeta)d\zeta\cdot \int_{0}^{x}\hspace{-0.1cm}z^{2}(\zeta,t)d\zeta\Big)dx\cr\hspace{-0.5cm}&&
\leq 2\mathbb{E}\int_0^1 z^2(x,t)dx+2\max_{0\leq x\leq 1}\max_{0\leq\zeta\leq x}l^{2}(x,\zeta)\cdot\mathbb{E}\int_0^1 z^2(x,t)dx
\cr\hspace{-0.5cm}&&= {\Gamma^{*}e^{-\theta^{*}t}},
\end{eqnarray}
where
\begin{equation}\label{constatsd}
\Gamma^{*}\triangleq 2(1+\max_{0\leq x\leq 1}\max_{0\leq\zeta\leq x}l^{2}(x,\zeta))\Gamma.
\end{equation}
\hfill$\Box$

\vspace{0.2cm}
{\bf APPENDIX E: Proof of Theorem \ref{fde67g}.}  Let $n=1,2,\cdots.$ Similarly  to the techniques in  (\ref{u-ito-dif}), (\ref{GWl2}),
it follows from It\^{o}'s formula  that for $n-1\leq t\leq n$,
\begin{eqnarray}
&& |z(\cdot,t)|^2_{L^2(0,1)} \cr &&=
|z(\cdot,n-1)|^2_{L^2(0,1)}+(\sigma^2-2c+2\varepsilon)\int^{t}_{n-1}|z(\cdot,s)|^2_{L^2(0,1)}ds
\cr &&+  \frac{1}{\varepsilon}\int^{t}_{n-1}\widetilde{w}^{2}(s)ds+\int^{t}_{n-1}2\sigma|z(\cdot,s)|_{L^2(0,1)}^2dB(s)\cr &&\leq
|z(\cdot,n-1)|^2_{L^2(0,1)}+\frac{1}{\varepsilon}\int^{t}_{n-1}\widetilde{w}^{2}(s)ds
+\int^{t}_{n-1}2\sigma|z(\cdot,s)|_{L^2(0,1)}^2dB(s),
\end{eqnarray}
where $0<\varepsilon<1$. Thus,
\begin{eqnarray}\label{fde3}
\hspace{-0.9cm}&& \mathbb{E}(\sup_{n-1\leq t\leq n}|z(\cdot,t)|^2_{L^2(0,1)})  \leq \mathbb{E}|z(\cdot,n-1)|^2_{L^2(0,1)}+
\frac{1}{\varepsilon}\int^{n}_{n-1}\mathbb{E}\widetilde{w}^{2}(s)ds\cr\hspace{-0.9cm}&&
+\mathbb{E}(\sup_{n-1\leq t\leq n}\int^{t}_{n-1}2\sigma|z(\cdot,s)|_{L^2(0,1)}^2dB(s)).
\end{eqnarray}

By the  Burkholder-Davis-Gundy inequality (see, e.g., \cite[Theorem 1.7.3, p. 40]{mao2007})
\begin{eqnarray}\label{bz1}
&&\mathbb{E}(\sup_{n-1\leq t\leq n}\int^{t}_{n-1}2\sigma|z(\cdot,s)|_{L^2(0,1)}^2dB(s))
\cr&& \leq 4\sqrt{2}\mathbb{E}\left(\int^{n}_{n-1}4\sigma^{2}|z(\cdot,s)|_{L^2(0,1)}^4ds\right)^{\frac{1}{2}} \cr&&
\leq 4\sqrt{2}\mathbb{E}\left(\sup_{n-1\leq s\leq n}|z(\cdot,s)|_{L^2(0,1)}^2\int^{n}_{n-1}4\sigma^{2}|z(\cdot,s)|_{L^2(0,1)}^2ds\right)^{\frac{1}{2}}
\cr&& \leq \frac{1}{2}\mathbb{E}(\hspace{-0.1cm}\sup_{n-1\leq t\leq n}|z(\cdot,t)|^2_{L^2(0,1)}\hspace{-0.05cm})\hspace{-0.05cm}+ 64\sigma^{2}\hspace{-0.1cm}\int^{n}_{n-1}\hspace{-0.2cm}\mathbb{E}|z(\cdot,s)|_{L^2(0,1)}^2\hspace{-0.06cm}ds.
\end{eqnarray}
After substitution of (\ref{bz1})  into (\ref{fde3}), we obtain  from Lemmas \ref{lem-Z-e-to0} and \ref{lemmas34} that
\begin{eqnarray}\label{dfe45}
\hspace{-0.7cm}&& \mathbb{E}(\sup_{n-1\leq t\leq n}|z(\cdot,t)|^2_{L^2(0,1)}) \leq 2\mathbb{E}|z(\cdot,n-1)|^2_{L^2(0,1)}\cr\hspace{-0.7cm}&&+
\frac{2}{\varepsilon}\int^{n}_{n-1}\mathbb{E}\widetilde{w}^{2}(s)ds+128\sigma^{2}\int^{n}_{n-1}\mathbb{E}|z(\cdot,s)|_{L^2(0,1)}^2ds
\cr\hspace{-0.7cm}&& \leq \Theta {e^{-\theta^{*}(n-1)}},
\end{eqnarray}
where $\Theta=(2+\frac{128\sigma^{2}}{\theta^{*}})\Gamma+\frac{2}{\varepsilon}|C|^{2}\frac{ \lambda^{2}_{\max}(Q_{c})}{\lambda_{\min}(Q_{c})(1-\mu\sigma^2)}\mathbb{E}[|Z(0)|^2+|\eta(0)|^2].$
Let $\varepsilon\in (0,\theta^{*})$ be arbitrary. By (\ref{dfe45}) and Chebyshev's inequality, it follows that
\begin{eqnarray}
\hspace{-1.3cm}&& P\{\sup_{n-1\leq t\leq n}|z(\cdot,t)|^2_{L^2(0,1)}> e^{-{(\theta^{*}-\varepsilon)}(n-1)}\} \cr\hspace{-1.3cm} && \leq
e^{{(\theta^{*}-\varepsilon)}(n-1)}\mathbb{E}(\sup_{n-1\leq t\leq n}|z(\cdot,t)|^2_{L^2(0,1)}) \leq
\Theta {e^{-\varepsilon(n-1)}}.
\end{eqnarray}

Applying   the Borel-Cantelli lemma (\cite[Lemma 2.4, p.7]{mao2007}), we obtain   for almost all $ \omega\in \Omega$, that
\begin{eqnarray}\label{bcdfe}
&&\sup_{n-1\leq t\leq n}|y(\cdot,t)|^2_{L^2(0,1)}\leq \Gamma^{*}\sup_{n-1\leq t\leq n}|z(\cdot,t)|^2_{L^2(0,1)}\cr&&\leq \Gamma^{*}e^{-{(\theta^{*}-\varepsilon)}(n-1)},
\end{eqnarray}
which holds for all but finitely many $n$ with $\Gamma^{*}$ given in (\ref{constatsd}).  Hence, there exists a random variable
 $n_{0}=n_{0}(\omega)$, such that for almost all $\omega\in \Omega$,
 (\ref{bcdfe}) holds whenever $n\geq n_{0}$. Hence, for almost all $\omega\in \Omega$,
\begin{eqnarray}
\frac{1}{t}\log|y(\cdot,t)|_{L^2(0,1)}=\frac{1}{2t}\log|y(\cdot,t)|^{2}_{L^2(0,1)}\cr \leq \frac{\log\Gamma^{*}}{2(n-1)}-\frac{{(\theta^{*}-\varepsilon)}(n-1)}{2n}
\end{eqnarray}
almost surely when $n-1\leq t\leq n$. Therefore,
\begin{equation}
\limsup_{t\to \infty}\frac{1}{t}\log|y(\cdot,t)|_{L^2(0,1)} \leq -\frac{{(\theta^{*}-\varepsilon)}}{2} \;\;\;\; \mbox{almost surely.}
\end{equation}
Since  {$\varepsilon>0$} is arbitrary, we then have
\begin{equation}
\limsup_{t\to \infty}\frac{1}{t}\log|y(\cdot,t)|_{L^2(0,1)} \leq -\frac{{\theta^{*}}}{2} \;\;\;\; \mbox{almost surely.}
\end{equation}
\hfill$\Box$

\vspace{-3mm}


\end{document}